\RequirePackage{fix-cm}
\documentclass[smallextended]{svjour3}       % onecolumn (second format)
\usepackage{amsmath}
\smartqed  % flush right qed marks, e.g. at end of proof
\usepackage{graphicx}
\usepackage{multirow}
\usepackage{epsfig}
\usepackage{epstopdf}
\usepackage{amsfonts}

\newcommand{\degree}{\ensuremath{^\circ}}
\journalname{}

\title{Triacontagonal proofs of the Bell-Kochen-Specker theorem}

\author{P.K.Aravind$^{1,}$$^{5}$, Justin Y.J.Burton$^{1}$,  Guillermo Núñez Ponasso$^{2,3}$ and D.Richter$^{4}$}

\institute{\at
$^{1}$Physics Department, Worcester Polytechnic Institute, Worcester, MA 01609, U.S.A.\\
$^{2}$Graduate School of Information Sciences, Division of Mathematics, Tohoku University, Sendai, Miyagi, Japan.\\
$^{3}$Department of Electrical and Computer Engineering, Worcester Polytechnic Institute, Worcester, MA 01609, U.S.A. \\
$^{4}$Department of Mathematics, Western Michigan University, Kalamazoo, MI 49008, U.S.A. \\
$^{5}$Corresponding author. \\
\email{paravind@wpi.edu},jyburton@wpi.edu,guillermo.carlo.nunez.a7@tohoku.ac.jp,gcnunez@wpi.edu,david.richter@wmich.edu}
\date{\today}

\begin{document}
\maketitle
\begin{abstract}
Coxeter pointed out that a number of polytopes can be projected orthogonally into two dimensions in such a way that their vertices lie on a number of concentric regular triacontagons (or 30-gons). Among them are the 600-cell and 120-cell in four dimensions and Gosset’s polytope $4_{21}$ in eight dimensions. We show how these projections can be modified into Kochen-Specker diagrams from which parity proofs of the Bell-Kochen-Specker theorem are easily extracted. Our construction trivially yields parity proofs of fifteen bases for all three polytopes and also allows many other proofs of the same type to be constructed for two of them. The defining feature of these proofs is that they have a fifteen-fold symmetry about the center of the Kochen-Specker diagram and thus involve both rays and bases that are multiples of fifteen. Any proof of this type can be written as a word made up of an odd number of distinct letters, each representing an orbit of fifteen bases. Knowing the word representing a proof makes it possible to infer all its characteristics without first having to recover its bases. A comparison is made with earlier approaches that have been used to obtain parity proofs in these polytopes, and some directions in which this work can be extended are discussed.\\ 

\end{abstract}
\section{\label{sec:Intro}Introduction}

Coxeter pointed out that a number of four and eight dimensional polytopes can be projected orthogonally into two dimensions in such a way that their vertices lie on a number of concentric regular triacontagons (or 30-gons). He dubbed this the triacontagonal projection and the frontispiece of his monograph \textit{Regular Polytopes}\cite{Cox1} shows such a projection of the 600-cell. Similar projections of the 120-cell and Gosset's polytope $4_{21}$ are shown in \cite{Cox2} and \cite{leonardo}, respectively. All three projections have $C_{30}$ point symmetry about their center, which gives them a highly symmetrical appearance. Table 1 lists the features of these polytopes that are of interest in connection with proofs of the Bell-Kochen-Specker (BKS) theorem\cite{BKS,KS}\footnote{This theorem rules out the existence of deterministic noncontextual hidden variables theories as an alternative to quantum mechanics. See Sec.2 for an overview of it.}, namely, the system of rays and bases they give rise to. A ray is an undirected line through the center of the polytope that passes through a pair of antipodal vertices, either of which may be taken as a representative of the ray. The number of rays of each of the polytopes, whose vertices come in antipodal pairs, is half the number of its vertices. A basis is a set of four or eight mutually orthogonal rays (four for the 600-cell and 120-cell and eight for Gossets's polytope). \\

\begin{table}[ht]
	\centering % used for centering table
	\begin{tabular}{|c | c | c|} % centered columns (2 columns)
		\hline % inserts single horizontal line
		\bf{Polytope} & \bf{Dimension} & \bf{Ray-Basis symbol} \\
		\hline 
		600-cell or $\{3,3,5\}$  & 4 & $60_{5} - 75_{4}$ \\
		\hline
		120-cell or $\{5,3,3\}$  & 4 & $300_{9} - 675_{4}$ \\
		\hline	
		Gosset's polytope or $4_{21}$ & 8 & $120_{135} - 2025_{8}$ \\
		\hline
	\end{tabular}
	\caption{The three polytopes studied in this paper, with their Schl\"{a}fli symbols (for the first two) or  Coxeter symbol (for the third). The dimension of the polytope is indicated in the second column and its ray-basis symbol in the third. The number to the left of the dash in the symbol is the number of rays of the polytope and the number to the right is the number of bases, with the subscript on the left indicating the number of bases to which each ray belongs and that on the right the number of rays in each basis.}
	\label{Poly} % is used to refer this table in the text
\end{table}

The purpose of this paper is to show how the triacontagonal projections of the three polytopes can be converted into Kochen-Specker diagrams from which parity proofs of the BKS theorem are easily extracted. It is interesting that these projections, which arose out of purely geometrical considerations, should lead so directly to an important physical result.  \\

The plan of this paper is as follows. Sec.2 gives a quick introduction to the BKS theorem and the notion of a parity proof that should allow readers unfamiliar with these topics to follow the rest of the (largely geometrical) arguments in this paper. Sec.3 introduces the triacontagonal projection and shows how it can be modified into a Kochen-Specker diagram from which parity proofs of the BKS theorem can be extracted. A Kochen-Specker diagram (or orthogonality graph) is a graph whose vertices are rays and whose edges connect orthogonal pairs of rays. An important feature of a Kochen-Specker diagram obtained from a triacontagonal projection is that both the rays \textbf{and the bases\footnote{This will become clear in Sec.3.}} have a fifteen-fold symmetry about the center of the diagram. As a result, it has parity proofs with a fifteen-fold symmetry about the center of the diagram \textbf{and it is only with such proofs that this paper will be concerned}. We will refer to these proofs as ``fifteen-fold proofs'' or simply proofs in the rest of this paper. There are a large number of other proofs of lower symmetry, or no symmetry at all, but we will largely ignore them except as they relate to the proofs discussed here.\\

All parity proofs in the three polytopes considered here have been found earlier by a variety of other methods \cite{Waegell2011b,Waegell2014,Waegell2015,Pavicic2011,Megill,Lisonek2014b}, with our fifteen-fold proofs included among them. However, the symmetry of the proofs went unrecognized and was not exploited in any way. We show that by recognizing the symmetry at the outset, it is possible to categorize the proofs better and also give a more detailed account of the members of particular subclasses. In particular, we show that the proofs of highest symmetry in all three polytopes can be written as ``words'' made up of a number of distinct letters from which all their properties can be inferred.\\

Sections 4,5 and 6 discuss the specific form the triacontagonal projection takes for each of the polytopes and how it facilitates the extraction of the fifteen-fold parity proofs in them. The task turns out to be trivial for the 600-cell, which has only a handful of such proofs, but is more involved in the other two cases in which the number of such proofs runs into the billions. We use an approach due to Lison{\v e}k et al\cite{lisonek1}, suitably modified to take advantage of the symmetry of the triacontagonal projection, to show how the proofs can be written as ``words'' made up of distinct letters, each of which represents an orbit of fifteen bases. We give a simple algorithm for constructing all the words representing arbitrary fifteen-fold parity proofs of the 120-cell and Gosset's polytope and show how all the characteristics of the proofs can be inferred from the words themselves without the need to recover the bases involved in the proof. Since many of the proofs consist of hundreds of bases, the benefits of such an approach are obvious.\\    

Section 7 recapitulates the main results of the paper and also indicates some directions in which this work might be extended.\\ 
  
Before going on, we give a brief overview of the field of quantum contextuality in which the present work is situated. The basic result that laid the foundations for the field was the discovery by Bell\cite{BKS} and Kochen and Specker\cite{KS} of the incompatibility between deterministic noncontextual hidden variables theories and quantum mechanics. Bell arrived at this conclusion by a continuum argument, but Kochen and Specker established it by using 117 rays in a three-dimensional Hilbert space. The latter approach has been greatly generalized over the years, and many examples are now known of rays in both real and complex Hilbert spaces of all dimensions that prove the theorem (see \cite{Budroni} for many references to such work). The three polytopes studied in this paper are just one of the many venues in which such proofs have been found. What makes them particularly interesting is that they each have well over a billion parity proofs in them, in addition to a vast number of other contextuality proofs that are not as easy to categorize or recognize.\\ 

Quantum contextuality has been shown to be useful for a variety of tasks such as ``all versus nothing" proofs\cite{Liu}, bipartite perfect quantum strategies\cite{BPQS}, the quantum computational advantage of shallow circuits\cite{Brayvi}, random number generation\cite{random} and universal quantum computation via magic state distillation\cite{Howard}. In a notable development, Spekkens \cite{Spekkens1} broadened the notion of a noncontextual theory by specifying the sorts of restrictions that would govern its predictions about the most general sorts of measurements that could be carried out on arbitrarily prepared systems. The advantage of such an approach is that it makes it possible to identify tasks in which a noncontextual theory predicts outcomes only within a certain range of parameters whereas experiment (and contextual models) allow values outside that range. This approach makes it possible to demonstrate the existence of quantum contextuality by means of robust experiments that are not subject to the sort of objections that had been made earlier to traditional contextuality tests\cite{Exp,Meyer,Clifton,Barrett}. An example of a practical  task in which a ``quantum advantage'' can be demonstrated by using Spekkens' approach is parity oblivious multiplexing\cite{Spekkens2}.\\

\section{\label{sec:BKS}The Bell-Kochen-Specker (BKS) theorem}

We give a quick introduction to the BKS theorem and the notion of a parity proof for readers unfamiliar with these topics so that they can follow the rest of the  paper (see \cite{Budroni} for a more detailed treatment of this topic and many references).\\

The BKS theorem rules out the existence of deterministic noncontextual hidden variables theories as an alternative to quantum mechanics. To understand  what such theories are, consider a set of rays in a Hilbert space of even dimension $N$ that form a number of bases, with each basis consisting of $N$ mutually orthogonal rays from the set. It is possible for a ray to occur in more than one basis of the set (and in the ray-basis sets of interest to us, which are shown in the third column of Table 1, this always happens). A noncontextual hidden variables theory is one that assigns the value 0 or 1 to every ray in the set, and further requires every ray to assume its assigned value in every basis in which it occurs. A contextual theory, on the other hand, allows a ray to switch its value between 0 and 1 in different bases. \\ 

It turns out that it is possible to do a laboratory experiment to determine whether each ray in a basis has the value 0 or 1, and that such an experiment always shows exactly one ray in any basis to have the value 1, with all the others having the value 0. In order to agree with this empirical finding, a noncontextual theory must be able to assign a 0 or 1 to every ray in a ray-basis set in such a way that every basis has just one ray assigned the value 1 in it. However, if there is a ray-basis set in an even number of dimensions with the properties that (a) the number of bases in the set is odd, and (b) each ray occurs in an even number of bases of the set, then the noncontextual theory cannot succeed in its task and must be rejected. The argument showing why it fails is this: every basis is required to have exactly a single 1 in it, and from (a) the total number of 1's over all the bases must be odd. However, from (b), the total number of 1's over all the bases must be even, and this contradicts the earlier conclusion. Because this argument relies on an even/odd contradiction, the ray-basis set it uses to achieve its goal is often referred to as a parity proof.\\

The purpose of the BKS theorem is to rule out the existence of noncontextual hidden variables theories, and one way of doing this is by finding a ray-basis set that satisfies the conditions stated above for a parity proof. That is not the only way of proving the theorem, but it is the method of choice for the polytopes considered in this paper.\\  

\section{\label{sec:Overview}From the triacontagonal projection to the Kochen-Specker diagram and the basis table}

The triacontagonal projection is an orthogonal projection of certain polytopes from their original dimension down to two dimensions in such a way that their vertices lie on a number of concentric regular triacontagons. Coxeter constructed such projections for the three polytopes studied here by obtaining triacontagonal coordinates for their vertices\footnote{The  triacontagonal coordinates of the 600-cell and 120-cell are given in \cite{Cox1} and those of Gosset's polytope in \cite{leonardo}.} with the  property that the first two coordinates specify the positions of the vertices in the plane of the projection and the remaining coordinates their positions in the subspace orthogonal to the projection plane.\\

The triacontagonal projection shows the vertices of the polytope as points and its edges as line segments joining the points. Two changes must be made to  convert it into a Kochen-Specker diagram (or orthogonality graph): firstly, every diametrically opposite pair of points in the projection should be replaced by just one of them  and, secondly, the line segments representing edges should be replaced by line segments joining orthogonal pairs of rays. The first change is needed because diametrically opposite triacontagon vertices in the projection correspond to antipodal vertices of the polytope, both of which represent the same ray. This change can be effected by keeping only the alternate vertices of all the triacontagons, thereby converting them into pentadecagons (or 15-gons).  The second change can be effected by using the   triacontagonal coordinates of the vertices to work out all the orthogonalities between the rays and then looking at just the first two coordinates to pick out all orthogonal pairs of rays in the projection plane. \\

The Kochen-Specker diagram  differs from its parent projection in having only half the number of points, a very different pattern of lines and a fifteen-fold symmetry about its center compared to the thirty-fold symmetry of the parent. However the visual appearance of the diagram is confusing and not particularly helpful, and we will make no direct use of it. Instead, we will replace it by a table of numbers, the so called \textbf{basis table}, that captures all the essential information in it. The basis table is simply a list of all the bases formed by the rays of the polytope. Since all the polytopes we are considering are saturated (i.e., have no orthogonalities between rays other than those contained in their bases), the basis table captures all the orthogonalities in the Kochen-Specker diagram and  so serves as a complete  substitute for it. However the basis table is more than just a convenient proxy for the Kochen-Specker diagram. As will become clear below, it facilitates the search for parity proofs, and also their viewing, in ways that the diagram simply cannot.\\

We will exploit the fifteen-fold symmetry of 
the Kochen-Specker diagram by introducing a numbering scheme for the rays that greatly simplifies the specification of the bases. We will number the rays in counterclockwise order around the pentadecagons, beginning from the outermost one and proceeding inwards towards smaller ones, always attaching the smallest number on any pentadecagon to the ray lying on, or just above, the positive x-axis. Then a rotation by one-fifteenth of a revolution (or $24^\circ$) about the center of the diagram sends any ray into the one whose number is one greater than it, with ray 1 going into ray 2, ray 24 into ray 25, and so on. However, a modification of this rule is clearly needed for rays that are multiples of fifteen, like 30 or 75, because such rays do not increase by one but go back into the rays at the beginning of their pentadecagons instead; thus ray 30  goes back into ray 16, ray 75 back into ray 51, and so on. We will term this  drop in the number of a ray as one goes around a pentadecagon as
``wraparound'' (see Table 4, and its caption, for a graphic illustration of this concept as one runs one's eye down the columns).\\

Next pick a set of rays in the diagram that form a basis\footnote{This is easily done using the triacontagonal coordinates of the rays.}. If one adds one to the numbers of each of the rays in this basis (with wraparound!), one can verify, using triacontagonal coordinates, that the new rays so obtained also form a basis. This  pattern holds no matter how often it is repeated and shows that the bases have a fifteen-fold symmetry about the center of the diagram, just like the rays. In view of this, it is possible to generate an orbit of fifteen bases by picking any one of its members and repeatedly adding one to each of its rays (with wrapround). We will term the representative picked from each orbit the \textbf{generator} and label it by a lower case letter.   The entire basis table can then be written as a ``word'' made up of a number of distinct letters, each representing a generator, and all the bases can be recovered by replacing each letter (or generator) by the orbit of fifteen bases it represents.\\ 

A fifteen-fold parity proof must consist of a number of complete orbits of bases (to respect its fifteen-fold symmetry), and so can be written as a word made up of the letters/generators representing its orbits of bases. Since the total   number of bases in a parity proof must be odd, and each letter contributes fifteen bases, the word representing a fifteen-fold proof must consist of an odd number of letters. However not all odd letter words are suitable because each ray is required to occur an even number of times over the bases. To see which odd letter words are suitable, it is convenient to introduce a capital letter for each ray in a basis to indicate the pentadecagon to which it belongs. Then every lower case letter in a word can be replaced by the set of capital letters representing the pentadecagons that cover its bases and the word can be written as a string of capital letters, each representing a complete pentadecagon's worth of fifteen rays. But the requirement that each ray occur an even number of times over the word implies that each capital letter also occur an even number of times. Thus the condition that an odd letter word must satisfy to be a parity proof can be stated as follows: replace each lower case letter in the word by its set of capital letters and count the number of times each capital letter occurs over the word; if all the counts are even the word represents a parity proof, otherwise not. A concrete illustration of this procedure is given in Eq.(1) for the 600-cell and Eq.(3) for the 120-cell, with the present explanation repeated there and supplemented by the details needed in the individual cases to make it clear.\\
  
In the next three sections we will describe an efficient method of picking out odd length words representing parity proofs in all three polytopes. \\     

We end with a word about terminology. Since the parity proofs of this paper have a fifteen-fold symmetry, they might best be termed ``pentadecagonal'' proofs. However we have referred to them as ``triacontagonal'' proofs in the title of this paper to acknowledge their origin in the projections of the same name due to Coxeter. 

\pagebreak

\section{\label{sec:600cell}The 600-cell}

The 600-cell is a four-dimensional regular polytope with 120 vertices lying on the surface of a 3-sphere. Its vertices come in antipodal pairs, so it has 60 rays. The rays form 75 bases\footnote{While the bases can be calculated using any of the standard  coordinates for the vertices, the advantage of using triacontagonal coordinates is that it allows the points corresponding to the bases to be picked out in the projection plane by looking at just the first two coordinates of the rays.}, with each ray occurring in five bases, so the ray-basis symbol of this polytope is $60_{5}-75_{4}$.\\

The triacontagonal projection of the 600-cell has 120 points representing its vertices and 720 line segments representing its edges. The Kochen-Specker diagram obtained from it has 60 points representing the rays and 450 line segments connecting orthogonal pairs of rays. The rays lie, in sets of fifteen, at the vertices of four concentric regular pentadecagons. Table 2 lists the numbering scheme we use for the rays.\\  

\begin{table}[ht]
	\centering % used for centering table
	\begin{tabular}{|c | c | } % centered columns (2 columns)
		\hline % inserts single horizontal line
		\bf{Pentadecagon} &  \bf{Rays}  \\
		\hline 
		A & 1-15 \\
		\hline
		B & 16-30 \\
		\hline
		C & 31-45 \\
		\hline
		D & 46-60 \\
		\hline
	\end{tabular}
	\caption{Pentadecagons of the 600-cell, labeled A through D in order of decreasing circumradius. The rays, which  lie at the vertices of the pentadecagons, are numbered in counter-clockwise order around them, with the lowest number on any  pentadecagon being attached to the ray lying on or just above the positive x-axis, and the numbers increasing as one moves inwards towards smaller pentadecagons. The pentadecagons A and D have vertices lying on the x-axis, whereas B and C have vertices lying  $6{\degree}$ above the x-axis.}
	\label{Tab2} % is used to refer this table in the text
\end{table}

The 600-cell has 75 bases that fall into five orbits of fifteen bases each. Table 3 gives a generator for each orbit along with its basis profile, with each capital letter of the profile indicating the pentadecagon to which the corresponding ray of the generator belongs.  \\ 
 
\begin{table}[ht]
	\centering % used for centering table
	\begin{tabular}{|c | c |} % centered columns (2 columns)
		\hline % inserts single horizontal line
		\bf{Basis profile} & \bf{Generator} \\
		\hline 
		AADD & $a = \{1,5,55,56\}$ \\
		\hline
		BBCC & $b = \{16,18,36,43\}$ \\
		\hline
		ABCD & $c = \{1,19,43,49\}$,  $d = \{1,20,41,58\}$, $e = \{1,27,42,46\}$ \\
		\hline
	\end{tabular}
	\caption{Generators of the 600-cell. The five generators are labeled $a$ through $e$ and the  four mutually orthogonal rays of which each is made up are indicated after it in curly brackets. The basis profile of each generator is indicated in the column to its left, with each capital letter indicating the pentadecagon to which the corresponding ray of the generator belongs. The first two generators have unique profiles, whereas the last three have a common profile. Despite this, the generators $c$,$d$ and $e$ give rise to orbits that have no bases in common (as is clear from Table 4).} 
	\label{Tab3} % is used to refer this table in the text
\end{table} 

\begin{table}[ht]
	\centering % used for centering table
	\begin{tabular}{|c c c c | c c c c | c c c c |  c c c c | c c c c |} % centered columns (5 columns)
		\hline % inserts single horizontal line
	\bf{A} & \bf{A} &\bf{D} &\bf{D} &\bf{B} &\bf{B} &\bf{C} &\bf{C} &\bf{A} &\bf{B} &\bf{C} &\bf{D} &\bf{A} &\bf{B} &\bf{C} &\bf{D} &\bf{A} &\bf{B} &\bf{C} &\bf{D}  \\
		\hline
		\bf{1} & \bf{5} &\bf{55} &\bf{56} &\bf{16} &\bf{18} &\bf{36} &\bf{43} &\bf{1} &\bf{19} &\bf{43} &\bf{49} &\bf{1} &\bf{20} &\bf{41} &\bf{58} &\bf{1} &\bf{27} &\bf{42} &\bf{46}  \\
		2 & 6 & 56 & 57 & 17 & 19 & 37 & 44 & 2 & 20 & 44 & 50 & 2 & 21 & 42 & 59 & 2 & 28 & 43 & 47  \\
		3 & 7 & 57 & 58 & 18 & 20 & 38 & 45 & 3 & 21 & 45 & 51 & 3 & 22 & 43 & 60 & 3 & 29 & 44 & 48  \\
		4 & 8 & 58 & 59 & 19 & 21 & 39 & 31 & 4 & 22 & 31 & 52 & 4 & 23 & 44 & 46 & 4 & 30 & 45 & 49  \\
		5 & 9 & 59 & 60 & 20 & 22 & 40 & 32 & 5 & 23 & 32 & 53 & 5 & 24 & 45 & 47 & 5 & 16 & 31 & 50  \\
		6 & 10 & 60 & 46 & 21 & 23 & 41 & 33 & 6 & 24 & 33 & 54 & 6 & 25 & 31 & 48 & 6 & 17 & 32 & 51  \\
		7 & 11 & 46 & 47 & 22 & 24 & 42 & 34 & 7 & 25 & 34 & 55 & 7 & 26 & 32 & 49 & 7 & 18 & 33 & 52  \\
		8 & 12 & 47 & 48 & 23 & 25 & 43 & 35 & 8 & 26 & 35 & 56 & 8 & 27 & 33 & 50 & 8 & 19 & 34 & 53  \\
		9 & 13 & 48 & 49 & 24 & 26 & 44 & 36 & 9 & 27 & 36 & 57 & 9 & 28 & 34 & 51 & 9 & 20 & 35 & 54  \\
		10 & 14 & 49 & 50 & 25 & 27 & 45 & 37 & 10 & 28 & 37 & 58 & 10 & 29 & 35 & 52 & 10 & 21 & 36 & 55  \\
		11 & 15 & 50 & 51 & 26 & 28 & 31 & 38 & 11 & 29 & 38 & 59 & 11 & 30 & 36 & 53 & 11 & 22 & 37 & 56  \\
		12 & 1 & 51 & 52 & 27 & 29 & 32 & 39 & 12 & 30 & 39 & 60 & 12 & 16 & 37 & 54 & 12 & 23 & 38 & 57  \\
		13 & 2 & 52 & 53 & 28 & 30 & 33 & 40 & 13 & 16 & 40 & 46 & 13 & 17 & 38 & 55 & 13 & 24 & 39 & 58  \\
		14 & 3 & 53 & 54 & 29 & 16 & 34 & 41 & 14 & 17 & 41 & 47 & 14 & 18 & 39 & 56 & 14 & 25 & 40 & 59  \\
		15 & 4 & 54 & 55 & 30 & 17 & 35 & 42 & 15 & 18 & 42 & 48 & 15 & 19 & 40 & 57 & 15 & 26 & 41 & 60  \\
	\hline
	\end{tabular}
	\caption{Basis Table of the 600-cell. The generators $a$ through $e$ are shown in bold font at the tops of the columns, just below their basis profiles (also in bold font), and the orbit of each generator is shown below it. Any basis in an orbit is obtained from the one above it by adding one to each of its rays, but with ``wraparound''. To understand what this means, consider the basis {25,27,45,37} in the second column. Adding one to each of its rays  gives {26,28,46,38}. However the third ray, 46, has overshot the boundary of the pentadecagon C to which its predecessor belonged, and so should be replaced by the ray 31 at the beginning of that pentadecagon. Several other instances of wraparound occur in the table. The 75 bases shown here make up the entire basis table of the 600-cell. The orbits of the generators $a$ and $b$ each give rise to a $30_{2}-15_{4}$ parity proof, but the orbits of the other generators do not.}     
	\label{Tab4} % is used to refer this table in the text
\end{table}  

Table 4 shows the orbits of fifteen bases arising from the generators $a$ through $e$ of Table 3. The 75 bases shown there make up the entire basis table of the 600-cell, which can be described by the word $abcde$. The bases associated with the generators $a$ and $b$ each yield a $30_2-15_4$ parity proof involving thirty rays that each occur twice over an odd number (fifteen) of bases.\\

One can ask if the 600-cell has any fifteen-fold proofs beyond the two found above. Any such proof must be expressible as a ``word'' made up of an odd number of distinct letters from $a,b,c,d$ and $e$ (to ensure that the total number of bases in it is odd). However, it must also have each of the pentadecagons A,B,C and D occur an even number of times over its letters to ensure that each ray occurs an even number of times in it. It is easy to see that these two conditions allow only the eight proofs described by the words $a,b, acd, ace, ade, bcd, bce$ and $bde$, of which the six three-letter words are the new ones. The ray-basis symbol associated with any three-letter proof can be obtained by replacing each of its letters/ generators by its basis profile, counting the number of times each of the pentadecagons occurs over the word and finally converting the pentadecagon counts into ray counts. We illustrate how this  works for the word $acd$ below:\\

$acd$ $\rightarrow$
(AADD)(ABCD)(ABCD) $\rightarrow$ (CD)$_{2}$(AB)$_{4}$ $\rightarrow$$30_{2}30_{4}-45_{4}$   \hspace{2mm}   (1)
\\

\noindent
All the fifteen-fold proofs of the 600-cell are shown in Table 6, along with their ray-basis symbols. 

\begin{table}[ht]
	\centering % used for centering table
	\begin{tabular}{|c | c | c |} % centered columns (2 columns)
		\hline % inserts single horizontal line
		\bf{Parity proof} & \textbf{Pentadecagon counts} & \bf{Ray-Basis symbol}\\
		\hline
		$a$ & (AD)$_{2}$ & 30$_{2}$--15$_{4}$ \\
		\hline
	   	$b$ & (BC)$_{2}$ & 30$_{2}$--15$_{4}$ \\
	    \hline
		$acd$, $ace$, $ade$ & (BC)$_{2}$(AD)$_{4}$ & $30_{2}30_{4}-45_{4}$ \\
		\hline
	    $bcd$, $bce$, $bde$ & (BC)$_{2}$(AD)$_{4}$ & $30_{2}30_{4}-45_{4}$ \\
	    \hline
		\end{tabular}
	\caption{Fifteen-fold parity proofs of the 600-cell. The first column shows the word(s) representing the proof(s), the second its/their pentadecagon counts and the third the ray-basis symbol. All proofs with the same entries in the second and third columns have been lumped together in the first column.}
	\label{Tab5} % is used to refer this table in the text
\end{table}

We will term a parity proof \textbf{irreducible} if no subset of its bases constitutes a parity proof of any kind (even one of less than fifteen-fold symmetry) and \textbf{reducible} otherwise. Of the proofs in Table 5, only the single letter proofs are irreducible while the six three-letter proofs are reducible because they can all be reduced to one of the single letter proofs $a$ or $b$ by dropping two of their letters. \\

\section{\label{sec:120cell}The 120-cell}

The 120-cell is a four-dimensional regular polytope with 600 vertices lying on the surface of a 3-sphere. Its vertices come in antipodal pairs, so it has 300 rays. The rays form 675 bases, with each ray occurring in nine bases, so the ray-basis symbol of this polytope is $300_{9}-675_{4}$. \\

The Kochen-Specker diagram obtained from the triacontagonal projection arranges the rays at the vertices of twenty concentric regular pentadecagons, sixteen of which lie in pairs on eight different circles and four of which lie on circles of their own. Following Chilton\cite{Chilton}, we will label the pentadecagons by the letters A through L in order of decreasing circumradius, but add subscripts to eight of the letters to distinguish between pentadecagons lying on circles of the same radius. The numbering scheme we use for the rays is shown in Table 6. The main difference from the 600-cell is that there are now pairs of pentadecagons on a common circumcircle, with the vertices of one interleaving  those of the other. In such cases, we let the numbers increase all the way around one pentadecagon before we let them increase around the other (see Table 6 for details). \\  

The 675 bases\footnote{If one calculates the bases using triacontagonal coordinates, then the generators in Table 7 can be picked out very easily.} of the 120-cell fall into forty five different orbits of fifteen bases each. Table 7 lists one generator for each orbit, along with its basis profile. Most generators have a unique profile, but there are two sets of three generators that each have a common profile.
All the generators give rise to distinct orbits, irrespective of whether they share a profile with others or not. \\

Four of the generators in Table 7, namely, the ones labeled $j,q,r'$ and $s'$ give rise to $30_2-15_4$ parity proofs.  However it is possible to obtain a large number of other proofs by combining the generators in Table 7 suitably, and we proceed to show  how this can be done. \\

\begin{table}[ht]
	\centering % used for centering table
	\begin{tabular}{|c | c | c| c|} % centered columns (2 columns)
		\hline % inserts single horizontal line
		\bf{Pentadecagon} &  \bf{Rays}  \\
		\hline 
		A  & 1-15 \\
		\hline
		B$_{1}$ & 16-30 \\
		\hline
		B$_{2}$ & 31-45 \\
		\hline
		C & 46-60 \\
		\hline
		D$_{1}$ & 61-75 \\
		\hline
		D$_{2}$ & 76-90 \\
		\hline
		E$_{1}$ & 91-105 \\
		\hline
		E$_{2}$ & 106-120 \\
		\hline
		F$_{1}$ & 121-135 \\
		\hline
		F$_{2}$ & 136-150 \\
		\hline
		G$_{1}$ & 151-165 \\
		\hline
		G$_{2}$ & 166-180 \\
		\hline
		H$_{1}$ & 181-195 \\
		\hline
		H$_{2}$ & 196-210 \\
		\hline
		I$_{1}$ & 211-225 \\
		\hline
		I$_{2}$ & 226-240 \\
		\hline
		J & 241-255 \\
		\hline
		K$_{1}$ & 256-270 \\
		\hline
		K$_{2}$ & 271-285 \\
		\hline
		L & 286-300 \\
		\hline
	\end{tabular}
	\caption{Pentadecagons of the 120-cell, labeled A through L in order of decreasing circumradius, with pairs having the same circumradius being distinguished from each other by a subscript. The rays are numbered in counter-clockwise order around the pentadecagons, beginning from the outermost one and proceeding inwards and with the lowest number on any pentadecagon being attached to the ray lying on or just above the positive x-axis. For pentadecagons lying on a common circumcircle, the numbers first increase around the one with the smaller subscript before they increase around the other.}
	\label{Tab6} % is used to refer this table in the text
\end{table}

\begin{table}[ht]
	\centering % used for centering table
	\begin{tabular}{|c|c||c|c|} % centered columns (4 columns)
		\hline % inserts single horizontal line
		\bf{Basis profile} & \textbf{Generator} & \textbf{Basis profile} &\textbf{Generator} \\
		\hline 
		AB$_{1}$K$_{2}$L & $a=\{1, 20, 274, 296\}$ & CD$_{2}$H$_{2}$K$_{2}$ & $z=\{46, 79, 201, 280\}$\\
		\hline
		AB$_{2}$K$_{1}$L & $b=\{1, 42, 268, 291\}$ & CD$_{2}$I$_{2}$J & $a'=\{46, 80, 235, 244\}$\\
		\hline
		ACJL & $c=\{1, 50, 245, 286\}$ & CE$_{1}$H$_{2}$J & $b'=\{46, 94, 200, 247\}$\\
		\hline
		AD$_{1}$I$_{1}$L & $d=\{1, 65, 222, 289\}$ & CE$_{2}$H$_{1}$J & $c'=\{46, 110, 184, 250\}$\\
		\hline
		AD$_{2}$I$_{2}$L & $e=\{1, 87, 230, 298\}$ & CF$_{1}$G$_{2}$J & $d'=\{46, 125, 168, 246\}$\\
		\hline
		AE$_{1}$H$_{2}$L & $f=\{1, 95, 207, 292\}$ & CF$_{2}$G$_{1}$J & $e'=\{46, 139, 156, 251\}$\\
		\hline
		AE$_{2}$H$_{1}$L & $g=\{1, 117, 185, 295\}$ & D$_{1}$D$_{2}$I$_{1}$I$_{2}$ & $f'=\{61, 80, 212, 229\}$\\
		\hline
		AG$_{1}$H$_{1}$I$_{1}$ & $h=\{1, 155, 192, 215\}$ & D$_{1}$E$_{1}$H$_{2}$I$_{1}$ & $g'=\{61, 94, 209, 215\}$\\
		\hline
		AG$_{2}$H$_{2}$I$_{2}$ & $i=\{1, 177, 200, 237\}$ & D$_{1}$E$_{2}$F$_{1}$L & $h'=\{61, 118, 131, 298\}$\\
		\hline
		B$_{1}$B$_{1}$K$_{1}$K$_{1}$ & $j=\{16, 20, 265, 266\}$ & D$_{1}$E$_{2}$H$_{1}$I$_{1}$ & $i'=\{61, 116, 193, 224\}$\\
		\hline
		B$_{1}$D$_{1}$I$_{2}$K$_{2}$ & $k=\{16, 65, 228, 275\}$ & D$_{1}$F$_{2}$G$_{1}$I$_{2}$ & $j'=\{61, 140, 162, 234\}$\\
		\hline
		B$_{1}$E$_{1}$H$_{1}$K$_{2}$ & $l=\{16, 102, 193, 280\}$ & D$_{2}$E$_{1}$F$_{2}$L & $k'=\{76, 94, 141, 289\}$\\
		\hline
		B$_{1}$E$_{2}$I$_{1}$J & $m=\{16, 117, 221, 246\}$ & D$_{2}$E$_{1}$H$_{2}$I$_{2}$ & $l'=\{76, 96, 199, 228\}$\\
		\hline
		B$_{1}$F$_{2}$G$_{2}$K$_{1}$ & $n=\{16, 139, 178, 259\}$, & D$_{2}$E$_{2}$H$_{1}$I$_{2}$ & $m'=\{76, 118, 183, 237\}$\\
		& $o=\{16, 140, 176, 268\}$, &  & \\
		& $p=\{16, 147, 177, 256\}$ &  & \\
		\hline
		B$_{2}$B$_{2}$K$_{2}$K$_{2}$ & $q=\{31, 35, 280, 281\}$ & D$_{2}$F$_{1}$G$_{2}$I$_{1}$ & $n'=\{76, 132, 170, 218\}$\\
		\hline
		B$_{2}$D$_{2}$I$_{1}$K$_{1}$ & $r=\{31, 87, 224, 267\}$ & E$_{1}$E$_{2}$H$_{1}$H$_{2}$ & $o'=\{91, 116, 187, 200\}$\\
		\hline
		B$_{2}$E$_{1}$I$_{2}$J & $s=\{31, 95, 231, 244\}$ & E$_{1}$F$_{1}$G$_{2}$H$_{1}$ & $p'=\{91, 124, 172, 192\}$\\
		\hline
		B$_{2}$E$_{2}$H$_{2}$K$_{1}$ & $t=\{31, 110, 199, 262\}$ & E$_{2}$F$_{2}$G$_{1}$H$_{2}$ & $q'=\{106, 148, 160, 200\}$\\
		\hline
		B$_{2}$F$_{1}$G$_{1}$K$_{2}$ & $u=\{31, 125, 155, 271\}$, & F$_{1}$F$_{1}$G$_{1}$G$_{1}$ & $r'=\{121, 123, 156, 163\}$\\
		& $v=\{31, 132, 156, 274\}$, &  & \\
		& $w=\{31, 133, 154, 283\}$ &  & \\
		\hline
		CD$_{1}$H$_{1}$K$_{1}$ & $x=\{46, 65, 183, 269\}$ & F$_{2}$F$_{2}$G$_{2}$G$_{2}$ & $s'=\{136, 138, 171, 178\}$\\
		\hline
		CD$_{1}$I$_{1}$J & $y=\{46, 64, 224, 253\}$ & & \\
		\hline
	\end{tabular}
\caption{Generators of the 120-cell. The generators, which are labeled by lower case letters, are shown in the second and fourth columns, with the basis of each shown in  curly brackets. The basis profile of each generator is shown to its left, with each capital letter indicating the pentadecagon to which the corresponding ray of the generator belongs. The generators $n$,$o$ and $p$ all have a common profile, as do $u$,$v$ and $w$. Thus there are just 41 different profiles for the 45 generators. Each generator gives rise to an orbit of fifteen bases and the union of all the orbits gives the 675 bases of the 120-cell.}
\label{Tab7} % is used to refer this table in the text
\end{table}

We use a technique due to Lison{\u e}k et al\cite{lisonek1}, suitably adapted to our needs, to solve the problem. Consider the matrix $M$ whose $i,j$-th element, $M_{i,j}$, is the number of times pentadecagon $i$ of Table 6 occurs in the basis profile of generator $j$ of Table 7 (we number the pentadecagons and generators in the order they occur in these tables)\footnote{It is here that our approach deviates from \cite{lisonek1}, which uses rays and bases in place of our pentadecagons and generators in constructing the matrix $M$. Thus our approach might be considered a coarse-grained version of their fine-grained one. While our approach is unable to capture all the parity proofs in the 120-cell, it compensates by offering a far more detailed view of the ones it does.}. The matrix $M$ has dimensions 20 x 45. A parity proof can be extracted from $M$ by picking an odd number of its columns in such a way that if they are stacked side by side, the sum of the numbers in each of the horizontal rows is even. The word describing the resulting proof has for its letters the generators corresponding to the columns that are picked. It remains only to explain how the right columns of $M$ are to be picked.\\

The problem can be solved by finding all 45-dimensional column vectors $X$, with elements of 0 or 1, for which the equation $MX = 0$ (mod 2) is satisfied, for then the non-vanishing elements of $X$ will indicate the columns of $M$ to be picked. The solutions to the preceding equation are all vectors in the nullspace of the matrix $M$ mod 2. We used Maple to find that the nullspace has dimension 30 and found a set of linearly independent vectors in it. Taking all possible linear combinations of these vectors, with coefficients of 0 or 1, gives $2^{30}$ solutions to the equation $MX = 0$ (mod 2). However only half these solutions (namely, the ones with an odd number of 1's as elements) give parity proofs. Thus the number of parity proofs is $2^{29}$, and each of them can be constructed as the union of the bases associated with the generators picked out by the non-vanishing elements of $X$.\\

\begin{table}[ht]
	\centering % used for centering table
	\begin{tabular}{|c | c |} % centered columns (2 columns)
		\hline % inserts single horizontal line
		1-letter words & $j$ \\
		\hline
		2-letter words & $jq$, $jr'$, $js'$,$no$, $np$, $uv$, $uw$ \\
		\hline
		4-letter words & $cdjy$, $ceja'$, $cfjb'$, $cgjc'$, $dejf'$, $dfjg'$, $dgji'$, $efjl'$, $egjm'$, $fgjo'$ \\
		\hline
		6-letter words & $abdekr$, $abfglt$, $bcklsx$  \\
		\hline
		8-letter words & $ahklmsuh'$  \\
		\hline
		10-letter words & $abhijlmnsq'$, $acdefgklmz$, $acfghimsud'$, $bcfghimnse'$, $fghijlmsup'$ \\
		\hline
		12-letter words & $abfghijkmnsj'$, $bdefgiklmnsk'$, $defghijkmsun'$ \\
		\hline
	\end{tabular}
	\caption{The 30 words representing a set of linearly independent vectors in the nullspace of the matrix $M$. The first word is of odd length, while all the others are of even length. The letters of any word are always written in alphabetical order, with the primed letters written after the unprimed ones.}
	\label{Tab8} % is used to refer this table in the text
\end{table}

\begin{table}[ht]
	\centering % used for centering table
	\begin{tabular}{|c | c | c |} % centered columns (2 columns)
		\hline % inserts single horizontal line
		\bf{Length} & \bf{Parity proof} & \bf{Ray-basis symbol}\\
		\hline
		1 & $j$, $q$, $r'$, $s'$ & $30_{2}$-$15_{4}$ (irreducible)\\
		\hline
		3 & $cdy$, $def'$, $efl'$, $fgo'$ & $90_{2}$-$45_{4}=30_{2}$-$15_{4}\oplus30_{2}$-$15_{4}\oplus30_{2}$-$15_{4}$\\
		\hline
		5 & $abkrf'$, $ablto'$ & $150_{2}$-$75_{4}=30_{2}$-$15_{4}\oplus  ...$ (5 times)  \\
		\hline
		7 & $abegkri'$ & $150_{2}30_{4}$-$105_{4}$ (irreducible)\\
		\hline
	    7 & $bdklsxy$ & $180_{2}15_{4}$-$105_{4}$
		(irreducible)\\
		\hline
		9 & $fghilmsup'$ & $180_{2}45_{4}$-$135_{4}$ 
		(irreducible)\\
		\hline
		9 & $abhilmnsq'$ & $210_{2}30_{4}$-$135_{4}=42_{2}6_{4}$-$27_{4}\oplus  ... $ (5 times) \\
	    \hline
		11 & $abfghikmnsj'$, $defghikmsun'$ & $195_{2}45_{4}15_{6}$-$165_{4}$
		(irreducible)\\
		\hline
		13 & $bcdefghiknszq'$ & $225_{2}30_{4}15_{6}15_{8}$-$195_{4}$
		(irreducible) \\
		\hline
		15 & $bcdelnszd'e'h'k'n'p'q'$ & $165_{2}120_{4}15_{6}$-$225_{4}=170_{2}10_{4}$-$95_{4}+ ... $ (3 times)\\
		\hline
	\end{tabular}
	\caption{Fifteen-fold proofs of the 120-cell of all word lengths from 1 to 15 obtained by combining words in Table 8 according to the rule in Eq.(2). Words/proofs with the same ray-basis symbol have been lumped together in the first column. If the proof is reducible, the way its ray-basis symbol splits up into those of the smaller proofs in it is indicated in the last column (see the text for  an explanation). The symbol $\oplus$ is used if the smaller proofs have no bases in common, whereas  + is used if they do (this happens only for the last proof).  See Table 10 for a blown-up version of the proof $cdy$, which shows how its bases split up into those of the three smaller proofs in it. Note: the dots in the last column indicate that the ray-basis symbol after the equals sign is to be added to itself the number of times indicated in brackets at the end, with the type of addition, $\oplus$ or +, being the same as before. }
	\label{Tab9} % is used to refer this table in the text
\end{table}

\begin{table}[ht]
	\centering % used for centering table
	\begin{tabular}{|c c c c | c c c c | c c c c |} 
		% centered columns (3 columns)
		\hline 
		% inserts single horizontal line
		\bf{A} & \bf{C} & \bf{J} & \bf{L} & \bf{A} & \bf{D}$_{1}$ & \bf{I}$_{1}$ & \bf{L} & \bf{C} & \bf{D}$_{1}$ & \bf{I}$_{1}$ & \bf{J} \\
		\hline
		\bf{1} & \bf{50} &\bf{245} &\bf{286} & \bf{1} &\bf{65} &\bf{222} &\bf{289} &{46} &{64} &{224} &{253}   \\
    	\it{2} & \it{51} &\it{246} &\it{287} &\it{2} &\textit{66} &\textit{223} &\it{290} &\bf{47} &\bf{65} &\bf{225} &\bf{254}   \\
		{3} & {52} &{247} &{288} &{3} &{67} &{224} &{291} &\it{48} &\it{66} &\it{211} &\it{255}   \\
		\bf{4} & \bf{53} &\bf{248} &\bf{289} &\bf{4} &\bf{68} &\bf{225} &\bf{292} &{49} &{67} &{212} &{241}   \\
		\it{5} & \it{54} &\it{249} &\it{290} &\it{5} &\it{69} &\it{211} &\it{293} &\bf{50} &\bf{68} &\bf{213} &\bf{242}   \\
		{6} & {55} &{250} &{291} &{6} &{70} &{212} &{294} &\it{51} &\it{69} &\it{214} &\it{243}   \\
		\bf{7} & \bf{56} &\bf{251} &\bf{292} &\bf{7} &\bf{71} &\bf{213} &\bf{295} &{52} &{70} &{215} &{244}   \\
		\it{8} & \it{57} &\it{252} &\it{293} &\it{8} &\it{72} &\it{214} &\it{296} &\bf{53} &\bf{71} &\bf{216} &\bf{245}   \\
		{9} & {58} &{253} &{294} &{9} &{73} &{215} &{297} &\it{54} &\it{72} &\it{217} &\it{246}   \\
		\bf{10} & \bf{59} &\bf{254} &\bf{295} &\bf{10} &\bf{74} &\bf{216} &\bf{298} &{55} &{73} &{218} &{247}   \\
		\it{11} & \it{60} &\it{255} &\it{296} &\it{11} &\it{75} &\it{217} &\it{299} &\bf{56} &\bf{74} &\bf{219} &\bf{248}   \\
		{12} & {46} &{241} &{297} &{12} &{61} &{218} &{300} &\it{57} &\it{75} &\it{220} &\it{249}   \\
		\bf{13} & \bf{47} &\bf{242} &\bf{298} &\bf{13} &\bf{62} &\bf{219} &\bf{286} &{58}  &{61} &{221} &{250}   \\
		\it{14} & \it{48} &\it{243} &\it{299} &\it{14} &\it{63} &\it{220} &\it{287} &\bf{59} &\bf{62} &\bf{222} &\bf{251}   \\
		{15} & {49} &{244} &{300} &{15} &{64} &{221} &{288} &\it{60} &\it{63} &\it{223} &\it{252}   \\
		\hline
	\end{tabular}
	\caption{The three-letter proof  $cdy$. The generators $c$,$d$ and $y$ are shown in bold font at the tops of the three columns, just below their basis profiles, and their orbits of fifteen bases are shown below them. The 45 bases in this table give a $90_{2}-45_{4}$ parity proof, but this proof is reducible and breaks up into the three smaller $30_{2}-15_{4}$ proofs shown in bold, italic and regular font. Each of the smaller proofs has a five-fold symmetry (in that adding three to the rays of any basis yields the next basis in the chain). The smaller proofs have no bases in common, so $cdy$ is written as a direct sum ($\oplus)$ of them in Table 9.}. 
	\label{Tab10} % is used to refer this table in the text
\end{table}   

The 30 linearly independent vectors in the nullspace of $M$ can each be written as a ``word'' made up of distinct letters representing the generators in Table 7. One choice for these words is shown in Table 8, with the words listed in order of increasing length and the letters of any word always written in alphabetical order, with the primed letters coming after the unprimed ones (this convention eliminates the confusion caused by encountering the same word written in different forms). Henceforth, by the ``length'' of a fifteen-fold proof, we will mean the number of letters in the word representing it. \\ 

We now give a simple algorithm for obtaining all the fifteen-fold parity proofs. It consists of the following two steps: (1) pick an arbitrary subset of the 29 even length words in Table 8 and combine them together in the manner explained in the next paragraph to obtain an even length word, and (2) combine this word with the letter $j$ to obtain an odd length word representing a parity proof. The bases of the proof can be obtained as the union of the orbits of the letters/generators in its word. \\

The rule for combining two words alluded to above is the following. Let $U$ and $V$ be the sets of letters in the two words. Then the word obtained by combining them is represented by the set $W$ that is the symmetric difference of the sets $U$ and $V$:  
\\

\hspace{20mm} $W = U \bigtriangleup V := (U \cup V) \setminus (U \cap V)$
\hspace{12mm}    (2)
\\

\noindent
Stated in words (pun unintended!), Eq.(2) says that the word obtained by combining two words is the union of their sets of letters, but with all letters common to the two words dropped. We will actually use this rule to combine any two words, whether they belong to Table 8 or not. \\

The first step of our algorithm can be executed with the aid of Eq.(2) as follows: one combines any two words from the chosen subset of Table 8 to obtain a new word and then combines it with one of the remaining words from the subset, and keeps proceeding in this way until all the words in the subset get used up. The new words that result at each step of this process all lie outside the list in Table 8, with the final (even length) word being the input for the second step of the algorithm.  \\

An example should help make the above algorithm clear. Suppose one picks just the two even length words $abdekr$ and $dgji'$ from Table 8. Then, in the first step of the algorithm, one combines these words to get the 8-letter word $abegjkri'$, and in the second step one combines it with $j$ to get the 7-letter word $abegkri'$ that represents a parity proof. The ray-basis symbol of this proof can be worked out as follows: \\

$abegkri'$ $\rightarrow$ (B$_{1}$B$_{2}$D$_{1}$D$_{2}$E$_{2}$H$_{1}$I$_{1}$I$_{2}$K$_{1}$K$_{2}$)$_{2}$(AL)$_{4}$ $\rightarrow$ $150_{2}30_{4}-105_{4}$   \hspace{6mm}    (3)
\\

\noindent
The transition from the word to its ray-basis symbol can be made in the same way as for the 600-cell in Eq.(1) by using the basis profiles in Table 7.  \\ 

The algorithm we have just given for obtaining fifteen-fold proofs can be rephrased in terms of adding vectors in GF($2^{45}$). To do this, we first write each of the vectors in the nullspace of $M$ as a string of 45 bits whose 1's represent the letters/generators present in the corresponding word. The operation of combining a subset of the words in Table 8 can then be carried out by adding the corresponding 45-bit strings bitwise modulo 2 and converting the resulting string back into the word it represents. If the second step of the algorithm is omitted, only even length words are obtained, but adding the second step converts the even words into an equal number of odd length words that span an affine subspace of the same dimension as the even words. \\   

Table 9 shows examples of proofs of word lengths one to fifteen obtained by combining the words in Table 8, with the last column indicating if the proof is irreducible or not. For each reducible proof, the  smaller proofs in it are indicated by their ray-basis symbols, with the $\oplus$ sign being used if the proofs have no bases in common and the + sign if they do. One can determine if a proof is irreducible or not by setting up a matrix $M$ like the one we used earlier, but between the individual rays and bases of the proof, and determining the dimension of its nullspace; if the dimension is 1 the proof is irreducible, whereas if it larger than 1 it is reducible and one can determine the bases of all the smaller proofs in it. Table 10 shows how the reducible proof $cdy$ of Table 9 splits up into three smaller proofs that have no bases in common (hence the $\oplus$ sign). Note also that each of the smaller proofs has a five-fold symmetry, as explained in more detail in the caption to the table.\\  

 The proofs in Table 9 are just a vanishingly small fraction of all the fifteen-fold proofs that exist. To get a better idea of the proofs as a whole, we note that because each proof can be written as a 45-bit string of 1's and 0's, the proofs are like the codewords of odd weight\footnote{The weight of a binary codeword is the number of 1's in it.} of a linear binary code $[45,30,1]$\footnote{An $[n,k,d]$ linear binary code consists of binary codewords of $n$ bits of which $k$ are used to transmit a message and the other $n-k$ are used to correct errors and $d$ is the minimum distance between codewords (the distance between two codewords is the number of places in which their bits differ).} spanned by the 30 codewords in Table 8. One can  use the MacWilliams identities \cite{MacW1,MacW2} to calculate the number of codewords of each odd weight (and hence the number of parity proofs of each odd word length) and finds that the number of proofs of word length 1,3,5,7 or 9 is 4,48,564,5116 or 42576, respectively. The number grows rapidly after that and peaks at well over a hundred million for proofs of word length 23 before falling off to zero after a length of 39 (see Table 1 of Supplementary Materials for a complete histogram of the codewords).  \\

How many of the above proofs are reducible? While it is possible to tell by a direct examination if any given proof is reducible or not, is there a general statement about reducibility that one can make about a large class of proofs without having to examine every one? It turns out there is. Recall that a reducible proof is one whose word can be made to pass into the word of a smaller proof by leaving out an even number of the letters in it (such proofs of the 600-cell were exhibited in Table 5). In coding theory, a binary codeword is said to be \textbf{minimal} if it cannot be made to pass into another codeword by replacing some of its 1's by 0's. The following result of Ashikhmin and Barg\cite{Barg} provides some information about which codewords might be minimal: \\ 

\noindent
\textbf{Proposition 1.} The weight $l$ of a minimal codeword in an $[n,k,d]$ linear binary code must satisfy the inequality $l \leq n-k+1$. \\ 

\noindent
For the code associated with the 120-cell, $n=45$ and $k=30$ and Proposition 1 implies that a minimal codeword must have weight $l \leq 16$. Translating this into the language of our proofs, we can conclude, at one stroke, that no fifteen-fold proof of word length 17 or larger can be irreducible. However the proposition is ambivalent about proofs of length 15 or less and only a detailed  examination will tell if any particular proof is irreducible or not. We have actually carried out such an examination of all the proofs in Table 9 and listed the results in the last column of that table.\\    

To sum up, while Proposition 1 implies that fifteen-fold proofs of greater than a certain length are necessarily reducible, it says nothing about proofs of shorter length.\\   

\section{\label{sec:Gosset's polytope }Gosset's polytope $4_{21}$ }

Gosset's polytope $4_{21}$ is a uniform polytope in eight dimensions with 240 vertices lying on the surface of a 7-sphere. Its vertices come in antipodal pairs, so it has 120 rays. The rays form 2025 bases, with each ray occurring in 135 bases, so the ray-basis symbol of this polytope is $120_{135}-2025_{8}$. \\

The Kochen-Specker diagram obtained from the triacontagonal projection arranges the rays at the vertices of eight concentric regular pentadecagons which we label by the letters A through H in order of decreasing circumradius. See Table 11 for the numbering scheme we use for the rays.\\ 

\begin{table}[ht]
	\centering % used for centering table
	\begin{tabular}{|c | c | c| c|}
		% centered columns (2 columns)
\hline
		% inserts single horizontal line
		\bf{Pentadecagon} &  \bf{Rays} \\ \hline
		                         A           & 1-15      \\ \hline
		                         B           & 16-30     \\ \hline
		                         C           & 31-45     \\ \hline
		                         D          & 46-60     \\ \hline
		                         E          & 61-75     \\ \hline
		                         F          & 76-90     \\ \hline
		                         G          & 91-105    \\ \hline
		                         H         & 106-120   \\ \hline
	\end{tabular}
\caption{Pentadecagons of Gosset's polytope, labeled A through H in order of decreasing circumradius. The rays are numbered in counterclockwise order around the pentadecagons, beginning from the ray lying on or just above the positive x-axis.}
\label{Tab11} % is used to refer this table in the text
\end{table}

The 2025 bases of Gosset's polytope break up into 135 orbits of fifteen bases each\footnote{If one numbers  the rays as in Table 11 and uses triacontagonal coordinates for them, it becomes easy to pick out the generators of the bases and also to visualize the bases in the projection plane.}, with the orbits having 33 different basis profiles. Table 2 of the Supplementary Materials lists a generator for each orbit, along with its basis profile. Since there are 135 generators but only 33 basis profiles, many of the profiles have more than one generator associated with them.\\

An examination of all the basis profiles shows that sixteen generators associated with eight of the profiles directly give rise to parity proofs (two of which are shown as the first two entries in Table 13). However, many other proofs can be constructed by putting the generators together in the right combinations. This can be done by a method similar to that we used for the 120-cell in Sec.4, but with the matrix $M$ now having dimensions 8 x 135. The matrix can be set up by using the information about the generators and their basis profiles given in Table 2 of the Supplementary Materials. A calculation shows that the nullspace of $M$ mod 2 has dimension 131 and a set of linearly independent vectors in it are shown as the words of one to five letters (representing  generators) in Table 3 of the Supplementary Materials. By combining these words together in all possible ways using the rule given in Eq.(2) and keeping only the odd length words, one can generate all the $2^{130}$ fifteen-fold parity proofs in this polytope. \\

We now exhibit a few of the proofs. To spare the reader the necessity of looking up the Supplementary Materials, we have listed a few of the generators in Table 12 and some of the proofs that can be constructed from them in Table 13. \\  

\begin{table}[ht]
	\centering % used for centering table
	\begin{tabular}{|c || c || c |} % centered columns (2 columns)
		\hline % inserts single horizontal line
		\bf{Generator} & \bf{Basis Profile} \\
		\hline
		$a_{1}=\{1,4,21,76,95,109,111,115\}$&AABFGHHH\\
		\hline
		$b_{1}=\{1,4,35,37,112,114,116,118\}$&AACCHHHH\\
		\hline
		$c_{1}=\{1,5,32,50,105,113,115,119\}$&AACDGHHH\\
		\hline
		$d_{1}=\{1,4,35,74,76,97,111,112\}$&AACEFGHH\\
		\hline
		$e_{1}=\{1,4,51,59,91,94,114,116\}$&AADDGGHH\\
		\hline
		$e_{2}=\{1,5,46,50,100,101,115,116\}$&AADDGGHH\\
		\hline
		$h_{1}=\{1,18,21,42,87,106,108,119\}$&ABBCFHHH\\
		\hline
		$i_{4}=\{1,20,27,46,86,101,106,115\}$&ABBDFGHH\\
		\hline
		$m_{1}=\{1,18,43,51,73,78,101,111\}$&ABCDEFGH\\
		\hline
		$m_{8}=\{1,20,41,51,70,76,103,108\}$&ABCDEFGH\\
		\hline
		$n_{5}=\{1,27,32,71,73,81,88,96\}$&ABCEEFFG\\
		\hline
		$c'_{1}=\{16,19,31,50,51,66,76,93\}$&BBCDDEFG\\
		\hline
		$e'_{2}=\{16, 40, 41, 43, 58, 61, 74, 81\}$&BCCCDEEF\\
		\hline
		\end{tabular}
	\caption{Generators of Gosset's polytope. The first column shows 12 of the 135 generators, along with their letter labels, with their basis profiles shown to their right. The subscripts on the letter labels are used to distinguish generators belonging to the same basis profile.}
	\label{Tab12} % is used to refer this table in the text
\end{table}

\begin{table}[ht]
	\centering % used for centering table
	\begin{tabular}{|c | c | c |} % centered columns (2 columns)
		\hline % inserts single horizontal line
		\bf{Parity proof} & \textbf{Pentadecagon counts} & \bf{Ray-Basis symbol}\\
		\hline
		$b_{1}$ & (AC)$_{2}$(H)$_{4}$ & 30$_{2}$15$_{4}$--15$_{8}$ \\
		\hline
		$ e_{1}$ & (ADGH)$_{2}$ & $60_{2}-15_{8}$ \\
    	\hline
	    $a_{1}c_{1}e'_{2}$ & (BDEFG)$_{2}$(AC)$_{4}$(H)$_{6}$ & $75_{2}30_{4}15_{6}-45_{8}$ \\
    	\hline
    	$a_{1}h_{1}n_{5}$ &(CEG)$_{2}$(ABF)$_{4}$(H)$_{6}$ & $45_{2}45_{4}15_{6}-45_{8}$ \\
    	\hline
    	$c_{1}h_{1}i_{4}$ & (CDFG)$_{2}$(AB)$_{4}$(H)$_{8}$  & $60_{2}30_{4}15_{8}-45_{8}$ \\
    	\hline
        $a_{1}c_{1}d_{1}h_{1}m_{1}$ & (DE)$_{2}$(BCFG)$_{4}$(A)$_{8}$(H)$_{12}$ & $30_{2}60_{4}15_{8}15_{12}-75_{8}$ \\
    	\hline
    	$a_{1}c_{1}h_{1}m_{8}c'_{1}$ & (E)$_{2}$(CDFG)$_{4}$(AB)$_{6}$(H)$_{10}$ & $15_{2}60_{4}30_{6}15_{10}-75_{8}$ \\
    	\hline
	\end{tabular}
	\caption{Fifteen-fold proofs of Gosset's polytope. The transition from the words/proofs in the first column to the ray-basis symbols in the third can be made as in Eqs.(1) and (3) by using the basis profiles given in Table 12. All the proofs shown here are irreducible (i.e., they do not have any smaller parity proofs in them.}
	\label{Tab13} % is used to refer this table in the text
\end{table}

\begin{table}[ht]
	\centering % used for centering table
	\begin{tabular}{|c c c c c c c c | c c c c c c c c |} 
		% centered columns (2 columns)
		\hline 
		% inserts single horizontal line
		\bf{A} & \bf{A} & \bf{D} & \bf{D} & \bf{G} & \bf{G} & \bf{H} & \bf{H} & \bf{A} & \bf{A} & \bf{D} & \bf{D} & \bf{G} & \bf{G} & \bf{H} & \bf{H} \\
		\hline
    	\bf{1} & \bf{4} & \bf{51} & \bf{59} & \bf{91} & \bf{94} & \bf{114} & \bf{116} & {1} &  {5} &  {46} &  {50} &  {100} & {101} & {115} & {116} \\     
    	\hline
    	\it{2} & \it{5} & \it{52} & \it{60} & \it{92} & \it{95} & \it{115} & \it{117} & {2} & {6} & {47} & {51} & {101} & {102} & {116} & {117} \\
    	\hline
        \bf{3} & \bf{6} & \bf{53} & \bf{46} & \bf{93} & \bf{96} & \bf{116} & \bf{118} & {3} & {7} & {48} & {52} & {102} & {103} & {117} & {118} \\
        \hline
        \it{4} & \it{7} & \it{54} & \it{47} & \it{94} & \it{97} & \it{117} & \it{119} & \bf{4} & \bf{8} & \bf{49} & \bf{53} & \bf{103} & \bf{104} & \bf{118} & \bf{119} \\
        \hline
        {5} & {8} & {55} & {48} & {95} & {98} & {118} & {120} & \it{5} & \it{9} & \it{50} & \it{54} & \it{104} & \it{105} & \it{119} & \it{120} \\
        \hline
        \bf{6} & \bf{9} & \bf{56} & \bf{49} & \bf{96} & \bf{99} & \bf{119} & \bf{106} & {6} & {10} & {51} & {55} & {105} & {91} & {120} & {106} \\
        \hline
        \it{7} & \it{10} & \it{57} & \it{50} & \it{97} & \it{100} & \it{120} & \it{107} & {7} & {11} & {52} & {56} & {91} & {92} & {106} & {107} \\
        \hline
        \bf{8} & \bf{11} & \bf{58} & \bf{51} & \bf{98} & \bf{101} & \bf{106} & \bf{108} & {8} & {12} & {53} & {57} & {92} & {93} & {107} & {108} \\
        \hline
        \it{9} & \it{12} & \it{59} & \it{52} & \it{99} & \it{102} & \it{107} & \it{109} & \bf{9} & \bf{13} & \bf{54} & \bf{58} & \bf{93} & \bf{94} & \bf{108} & \bf{109} \\
        \hline
        {10} & {13} & {60} & {53} & {100} & {103} & {108} & {110} & \it{10} & \it{14} & \it{55} & \it{59} & \it{94} & \it{95} & \it{109} & \it{110} \\
        \hline
        \bf{11} & \bf{14} & \bf{46} & \bf{54} & \bf{101} & \bf{104} & \bf{109} & \bf{111} & {11} & {15} & {56} & {60} & {95} & {96} & {110} & {111} \\
        \hline
        \it{12} & \it{15} & \it{47} & \it{55} & \it{102} & \it{105} & \it{110} & \it{112} & {12} & {1} & {57} & {46} & {96} & {97} & {111} & {112}\\
        \hline
        \bf{13} & \bf{1} & \bf{48} & \bf{56} & \bf{103} & \bf{91} & \bf{111} & \bf{113} & {13} & {2} & {58} & {47} & {97} & {98} & {112} & {113} \\
        \hline
        \it{14} & \it{2} & \it{49} & \it{57} & \it{104} & \it{92} & \it{112} & \it{114} & \bf{14} & \bf{3} & \bf{59} & \bf{48} & \bf{98} & \bf{99} & \bf{113} & \bf{114} \\
        \hline
        {15} & {3} & {50} & {58} & {105} & {93} & {113} & {115} & \it{15} & \it{4} & \it{60} & \it{49} &\it{99} & \it{100} & \it{114} & \it{115} \\
        \hline
   	\end{tabular}
	\caption{Parity proofs in the word $e_{1}e_{2}$. The generators $e_{1}$ and $e_{2}$ are shown at the top of the left and right columns and their orbits of bases are shown below them. The bases in the left and right blocks each give a $30_{2}-15_{4}$ parity proof. However, there are three smaller parity proofs that get their bases from both columns of the table. One of them is shown in bold font and another in italic font. These two proofs have no bases in common and they both have the ray-basis symbol $36_2-9_{8}$. The third proof, which also has the same symbol, has some bases in common with the other two proofs but also three bases that are its own. If  the bases are numbered 1 to 30 from up to down and left to right, the bases associated with the third proof are 1,4,6,9,11,14,17,22 and 27.}. 
	\label{Tab14} % is used to refer this table in the text
\end{table}   

Table 14 shows all the bases associated with the two letter word  $e_{1}e_{2}$, which is not a parity proof but very interesting all the same. The bases associated with the letters $e_{1}$ and $e_{2}$ are shown in the two halves of the table, and the fifteen bases of each give a parity proof. However there are three smaller parity proofs in the table, each of which inherits its bases both from $e_{1}$ and $e_{2}$ (see table caption for details). All three proofs consist of just nine bases each and are the smallest parity proofs in Gosset's polytope. This example, like some of those in Table 9, illustrates the  fact that words of both  even and odd lengths can sometimes have much smaller parity proofs (of a less than fifteen-fold symmetry) embedded in them.\\ 

To round out the very limited view we have given of fifteen-fold proofs in Gosset's polytope, we used the MacWilliams identities to calculate the number of such proofs of various word lengths and found that the number of proofs of word length 1,3 and 5 is 16, 25812 and a little over 21 million, respectively, with the number growing rapidly after that and then falling off and going to zero beyond a word length of 133 (see Table 4 of the Supplementary Materials for more details). The proofs have a great variety of taxonomies, as reflected in their ray-basis symbols.\\

We finally address the question of which of all the  fifteen-fold proofs might be reducible. We can answer this question by using the fact that the fifteen-fold proofs of Gosset's polytope are like the codewords of an $[n,k,d]$ linear binary code with $n = 135$ and $k = 131$. Then it follows from Proposition 1 that all proofs with a word length of 7 or larger are reducible and that only the proofs with a word length of five or less might be irreducible. As we pointed out above, there are a little over 21 million proofs of five letters and we know of no way, short of an individual examination of each case, of telling how many of these might be irreducible.  \\ 

\section{\label{sec:Disc} Discussion}

We have shown how the triacontagonal projections of the 600-cell, 120-cell and Gosset's polytope can be used to obtain parity proofs of the BKS theorem having a fifteen-fold symmetry about the center of a suitably chosen Kochen-Specker diagram. The smallest such proofs in all three polytopes, which consist of fifteen bases, follow immediately. However there are many larger proofs that all involve both rays and bases that are multiples of fifteen.  We have shown that all such proofs can be written as words made up of an odd number of distinct letters from which all their properties, including the bases of which they are made up, can be inferred. Since the proofs number in the billions and many of them consist of hunderds of bases, the benefits of such an approach are obvious. While our construction captures only the proofs of highest symmetry in all three polytopes, it allows them to be categorized and accessed in a way that was not possible before.\\

The proofs studied in this paper all have a $C_{15}$ symmetry, which they inherit from the $C_{30}$ symmetry of the triacontagonal projection. The automorphism groups of all the polytopes are actually much larger\footnote{The symmetry groups of the 600-cell and 120-cell are identical, while that of Gosset's polytope is different. The symmetry groups of all three polytopes are discussed in \cite{Cox1}.}, but none of them has any rotational symmetries larger than $C_{30}$. However, since  $C_{5}$,  $C_{3}$ and trivial symmetry are all subgroups of $C_{15}$, these polytopes have many parity proofs with these smaller symmetries as well.\\

Aside from the three polytopes studied in this paper, we are aware of a family of eight-dimensional polytopes, all closely related to Gosset's polytope, for which triacontagonal projections exist. These polytopes can all be described by the Coxeter diagram below
\begin{figure}[htp]
	\begin{center}
		\includegraphics[width=.60\textwidth]{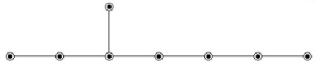}
	\end{center}
	\end{figure}
\noindent
by replacing one or more of its nodes by ringed nodes in all possible ways. Since this can be done in $2^{8}-1 = 255$ ways, there are 255 polytopes that can be obtained in this way, of which Gosset's polytope is one. However the remaining 254 polytopes also possess triacontagonal projections and any bases formed by their rays will therefore occur in orbits of fifteen. Thus it is possible that they also harbor parity proofs similar to those we have found in Gosset's polytope. However, working out the triacontagonal coordinates of these polytopes and seeing if this expectation is satisfied requires a lot of computation, and we leave it as a problem for the future.  \\

The fifteen-fold symmetry of the triacontagonal projection that we exploited in this paper suggests that a similar approach could prove useful in other cases. We mention two problems that could be tackled with this approach.  First, consider the system of 60 rays and 105 bases arising from the  triads of commuting observables of the two-qubit Pauli group. This system has been shown to contain a large number of parity proofs \cite{Waegell2011,lisonek1}, but only a few of them have been exhibited and no trends underlying them have been detected. It is possible that by looking at the automorphism group of this system and its subgroups, one could classify the proofs into different symmetry classes and achieve a much better understanding of them, similar to what has been done here. A far more challenging problem would be to apply this approach to the binary and ternary Golay codes, both of which are known to contain contextuality proofs that are a variant of parity proofs\cite{Waegell2022}. However, only one such proof has been identified in each case and neither is particularly small. The large dimensionalities of these systems (24 for the binary code and 12 for the ternary one) and the huge number of bases involved make the usual types of search methods impractical. However it is possible that by looking at the symmetry groups of these systems and identifying the right subgroups, one might be able to uncover at least a subclass of the contextuality proofs in each of them. Just as the triacontagonal projection yielded the smallest fifteen-fold proofs in all three polytopes trivially, it is our hope that symmetry could lead to a similar benefit in these far more intractable cases.   \\

We end by mentioning two open questions about the polytopes studied in this paper. The first concerns the size of the largest noncontextual sets in them, i.e., the largest number of bases with the property that every ray in them can be assigned a 0 or a 1 in such a way that every basis has exactly one ray assigned the value 1 in it. This question is interesting because the gap between the number of bases in the largest noncontextual set and the total number of bases of the polytope provides a measure of how contextual (or nonclassical) the bases of the polytope are. It was shown in \cite{BPQS} that there is a two-way connection between a quantum contextual set and a perfect quantum strategy in a nonlocal two player game based on that set. Determining the largest noncontextual set in a polytope would convey an idea of how large a quantum advantage it would provide in such a game. Unfortunately, we know of no good (i.e., efficient) algorithm for solving this problem and would welcome any progress  that could be made in this direction.\\ 

The second question concerns rigidity. The rays (or the equivalent projectors) of a polytope are said to be rigid if any set of projectors (not necessarily of rank one and possibly in a dimension greater than that of the polytope) that satisfy the same orthogonality relations as them can be mapped into them in a one-to-one fashion by a suitable unitary transformation. Proving the rigidity of a set of rays is not an easy task and not many examples of such sets are known. Rigidity is of interest because only a KS set that is rigid can be Bell self-tested\cite{Mayers}-\cite{Xu}. The discovery of two new rigid KS sets in $\mathbb{C}^{3}$ was recently reported in  \cite{Trandafir}, with one of them (consisting of 31 vectors) conjectured to be the smallest in three dimensions.\\

We would like to raise the question of whether the rays of the three polytopes
studied here are rigid or not. First consider the 600-cell. Each ray 
is an element of $\mathbb{RP}^{3}$ and  determined by three real 
numbers, and all 60 rays by 180 real numbers. However, the 75 bases 
of the 600-cell imply 450 orthogonality relations between its rays, and since 
the number of constraints imposed by orthogonality (450) greatly exceeds the number 
of ray components (180), it would seem that the rays are uniquely determined by orthogonality. 
However, there could be rays in $\mathbb{RP}^{n}$ or $\mathbb{CP}^{n}$ for 
$n > 3$ that obey the same orthogonality relations as the rays but are not unitarily equivalent to them, and this would demonstrate that the
600-cell is not rigid. We show in the Appendix that a subset of the rays of   Gosset's polytope satisfy all the orthogonality relations of the rays of the 
600-cell, but are not unitarily equivalent to them, and this establishes that the 
600-cell is not rigid.\\

However, we cannot make a statement about the rigidity of either the 120-cell or Gosset's polytope. The number of orthogonalities in both these polytopes (particularly the latter) is so large and imposes so many constraints that we feel that rigidity might well be forced. However, this is just a conjecture on our part and we have not come close to proving it. We leave it as an open problem to be considered. If either of these polytopes could be shown to be rigid, that would certainly enhance its interest as far as quantum contextuality is concerned.\\

One final point should be made. The bases of all three polytopes can be obtained as the products of powers of three unitary operators acting on the computational basis\cite{Waegell2011b}-\cite{Waegell2015}. If gates corresponding to the basic operators can be constructed, it would be possible to realize these polytopes (i.e., their bases) experimentally in the laboratory. This might be worth keeping in mind for the future time when it becomes possible. \\

{\bf Appendix: Proof that the 600-cell is not rigid}\\ 

We follow Moody and Patera \cite{Moody} in obtaining the 240 root vectors of $E_8$ 
(or the vertices of Gosset's polytope) using icosians, which are quantities related 
to the symmetries of the regular icosahedron.  We first construct the 120 root vectors 
of $H_4$ (which comprise the vertices of a regular 600-cell) in two different ways. 
Let $\beta=(1+\sqrt{5})/2$ be the golden ratio and $\alpha=(1-\sqrt{5})/2$ its conjugate.
Next, let $\mathbb{Z}[\alpha]=\{m+n\alpha: m,n \in \mathbb{Z}\}$ denote the ``golden ring'' 
consisting of all linear combinations of 1 and $\alpha$ with integer coefficients. 
In particular, we notice there is a bijection $\phi:\mathbb{Z}[\alpha]\rightarrow\mathbb{Z}^2$
which maps $m+n\alpha$ to $(m,n)$ and that this bijection extends to a map 
$\phi:\mathbb{Z}[\alpha]^4\rightarrow\mathbb{Z}^8$. \\

Let $G \subset GL(4,\mathbb{R})$ be the linear group generated by all even permutations 
and arbitrary sign changes of the four coordinates of a vector $(w, x, y, z) \in \mathbb{R}^{4}$. 
This group has $|G|= 12 \cdot 16 = 192$ elements. The roots of $H_4$ can be constructed as the union of the orbits of $(2,0,0,0)$, $(1,1,1,1)$ and $(0,\alpha,1,\beta)=(0,\alpha,1,-\alpha+1)$ under the 
action of $G$. The 120 vectors thus obtained represent the vertices of a 600-cell on a 
sphere of radius 2. The vectors come in pairs that are the negatives of each other, and 
if we keep just one member of each pair we get the 60 rays of the 600-cell. Let us call 
this set of 60 vectors (or rays) $H_4^a$. If we multiply each of the vectors of $H_4^a$ by 
$\alpha$ we get a set $H_4^b$ which serves as the set of vertices of a second 600-cell,
scaled by a factor of $\alpha$ relative to the first.  Analogous to $H_4^a$, the vectors of $H_4^b$ can be obtained as the union of the orbits of $(2\alpha,0,0,0)$, $(\alpha,\alpha,\alpha,\alpha)$ and $(0,1+\alpha,\alpha,1)$. \\
% I edited this paragraph heavily on September 10, 2025.  -DR

We define the 240 root vectors of $E_8$ as $\phi(H_4^a\cup H_4^b)$.
In more detail, we recall that every vector of either $H_4^a$ or $H_4^b$ has 
components that are elements of the golden ring $\mathbb{Z}[\alpha]$ and so can be written as 
$v=(m_{1}+n_{1}\alpha,m_{2}+n_{2}\alpha,m_{3}+n_{3}\alpha,m_{4}+n_{4}\alpha)$. 
Using the map $\phi$ defined above, we map this to the 8-dimensional vector with components 
$\phi(v)=(m_{1},m_{2},m_{3},m_{4},n_{1},n_{2},n_{3},n_{4})$. The 60 vectors obtained in this 
way from $H_4^a$ along with the 60 vectors obtained from $H_4^b$, together with all 
their negatives, make up the 240 root vectors of $E_8$ (and also the vertices of 
Gosset's polytope), cf.\ \cite{Moody}. \\
% I edited this paragraph heavily on September 10, 2025.  -DR

The 60 vectors of $H_4^a$ are unitarily equivalent to the 60 rays of the 600-cell 
introduced in Sec.3 and thus satisfy all the orthogonalites represented by the 75 bases 
of the 600-cell. If we take four vectors of $H_4^a$ that are mutually orthogonal and 
turn them into vectors of $E_8$ by the construction above, we find that the vectors 
we get in $E_8$ are also mutually orthogonal (in an 8-dimensional space.) However, the 
converse of this statement is not true, i.e., if one takes four vectors of $E_8$ that
are mutually orthogonal and looks at the vectors of $H_4^a$ (or $H_4^b$) that 
correspond to them, we will find that these latter vectors are generally not mutually 
orthogonal. \\

Here is an example that illustrates this. Consider the following four 
vectors of $H_4^a$: \\ 

\noindent
$v_1 = (2,0,0,0)$, \hspace{3mm}  $v_2=(0,\alpha,1,\beta)$,   
\hspace{3mm} $v_3=(0,1,\beta,\alpha,)$,  \hspace{3mm} $v_4=(0,\beta,\alpha,1)$    \\

\noindent
The four vectors of $E_8$ that correspond to these are \\

\noindent
$\phi(v_1) = (2,0,0,0,0,0,0,0)$, \hspace{1mm}  $\phi(v_2)=(0,0,1,1,0,1,0,-1)$,   \\
\hspace{1mm} $\phi(v_3)=(0,1,1,0,0,0,-1,1)$,  \hspace{1mm} $\phi(v_4)=(0,1,0,1,0,-1,1,0)$.    \\

\noindent
It can be checked that the vectors $\{v_1,v_2,v_3,v_4\}$ are mutually orthogonal, 
as are their images $\{\phi(v_1),\phi(v_2),\phi(v_3),\phi(v_4)\}\subset E_8$. Now consider the 
two additional vectors $v_5=(0,\alpha,1,-\beta)$ and  $v_6=(\alpha,\beta,1,0)$ of $H_4^a$ whose 
images in $E_8$ are $\phi(v_5)= (0,0,1,-1,0,1,0,1)$ and $\phi(v_6)=(0,1,1,0,1,-1,0,0)$.  One may
check that while the vectors $\{\phi(v_1),\phi(v_2),\phi(v_5),\phi(v_6)\}\subset E_8$
are mutually orthogonal,  their preimages $\{v_1,v_2,v_5,v_6\}\subset H_4^a$
are not because the vectors $v_1$ and $v_6$ are not orthogonal (and neither are $v_2$ and $v_5)$. \\
% I edited this paragraph heavily on September 12, 2025, including a correction to the vector 5'.  -DR

To summarize, we have constructed a bijective mapping between the rays of $E_8$ and 
$H_4^a\cup H_4^b$ with the property that all the orthogonalities between rays in $H_4^a$ 
are maintained between the corresponding rays of $E_8$, but with some orthogonalities 
between rays in $E_8$ not being preserved as orthogonalities between the corresponding 
rays of $H_4^a$. This demonstrates that the mapping we have constructed between $H_4^a$ 
and $E_8$ is not unitary (or orthogonal) and, therefore, that the 600-cell is not rigid.\\   

\noindent
\textbf{Acknowledgements.} 
We would like to thank both the anonymous referees for the many constructive suggestions they made for improving the manuscript.
GNP was supported by the Japanese Society for the Promotion of Science (JSPS) as a JSPS post-doctoral fellow and by the JSPS KAKENHI Grant number 24KF0176.

%\bibliographystyle{ieeetr}
%\bibliography{testbib}

\end{document}

\maketitle

\begin{table}[ht]
	\centering % used for centering table
	\begin{tabular}{|c | c | c| c |} % centered columns (2 columns)
		\hline % inserts single horizontal line
		\bf{Length} & \bf{Proofs} & \bf{Length} & \bf{Proofs }\\
                \hline
                1 &  4 & 21 &  115880440 \\
                \hline
                3 &  48 & 23 &  127058600 \\
                \hline
                5 &  564 & 25 &  96649232 \\
                \hline
                7 &  5116 & 27 &  51402240 \\
                \hline
                9 &  42576 & 29 &  19478640 \\
                \hline
                11 &  363392 & 31 &  5267856 \\
                \hline
                13 &  2371056 & 33 &  975268 \\
                \hline
                15 &  10593552 & 35 &  127856 \\
                \hline
                17 &  33022968 & 37 &  14708 \\
                \hline
                19 &  73615584 & 39 &  1212 \\
                \hline
        \end{tabular}
	\caption{Fifteen-fold symmetric parity proofs of the 120-cell. Shown above are the number of proofs for each word of odd length from 1 to 39. The total number of odd length proofs is $2^{29}$.}
	\label{Tab11} % is used to refer this table in the text
\end{table}

\begin{table}[ht]
	\centering % used for centering table
	\begin{tabular}{|c | c |} % centered columns (2 columns)
		\hline % inserts single horizontal line
		\bf{Basis profile} & \bf{Generator(s)} \\
		\hline 
		
		AABFGHHH & $a_{1}=\{1, 4, 21, 76, 95, 109, 111, 115\}$, $a_{2}=\{1, 4, 21, 81, 100, 114,\
		
		115, 116\}$, $a_{3}=\{1, 4, 21, 86, 105, 106, 115, 119\}$\\
		\hline
		AACCHHHH & $b_{1}=\{1, 4, 35, 37, 112, 114, 116, 118\}$\\
		\hline
		AACDGHHH & $c_{1}=\{1, 5, 32, 50, 105, 113, 115, 119\}$, $c_{2}=\{1, 5, 41, 46, 96, 112,\
		
		116, 118\}$\\
		\hline
		AACEFGHH & $d_{1}=\{1, 4, 35, 74, 76, 97, 111, 112\}$, $d_{2}=\{1, 4, 35, 75, 88, 105,\
		
		106, 112\}$, $d_{3}=\{1, 4, 37, 72, 89, 95, 109, 118\}$, \\
		& $d_{4}=\{1, 4, 37, 73, 86, 103, 118, 119\}$,\
		
		$d_{5}=\{1, 5, 32, 73, 85, 96, 118, 119\}$, $d_{6}=\{1, 5, 32, 74, 87, 94, 108, 119\}$,\
		
		\\
		& $d_{7}=\{1, 5, 41, 74, 76, 92, 108, 112\}$,\
		
		$d_{8}=\{1, 5, 41, 75, 78, 105, 112, 113\}$\\
		\hline
		AADDGGHH & $e_{1}=\{1, 4, 51, 59, 91, 94, 114, 116\}$, $e_{2}=\{1, 5, 46, 50, 100, 101,\
		
		115, 116\}$\\
		\hline
		AADEFGGH & $f_{1}=\{1, 4, 51, 72, 89, 91, 105, 106\}$, $f_{2}=\{1, 4, 51, 73, 76, 91, 103,\
		
		111\}$, $f_{3}=\{1, 4, 59, 74, 86, 94, 97, 119\}$, \\
		& $f_{4}=\{1, 4, 59, 75, 88, 94, 95, 109\}$,\
		
		$f_{5}=\{1, 5, 46, 64, 87, 91, 94, 116\}$, $f_{6}=\{1, 5, 48, 64, 89, 91, 105, 113\}$, \\
		& $f_{7}=\{1, 5, 48, 69, 89, 95, 96, 118\}$,\
		
		$f_{8}=\{1, 5, 48, 74, 89, 100, 101, 108\}$, $f_{9}=\{1, 5, 50, 69, 76, 92, 95, 115\}$\\
		\hline
		AAEEFFGG & $g_{1}=\{1, 4, 72, 74, 81, 89, 97, 100\}$, $g_{2}=\{1, 4, 73, 75, 81, 88, 100,\
		
		103\}$, $g_{3}=\{1, 5, 64, 73, 76, 85, 91, 92\}$, \\
		& $g_{4}=\{1, 5, 69, 75, 78, 87, 94, 95\}$,\
		
		$g_{5}=\{1, 5, 73, 75, 78, 85, 100, 101\}$\\
		\hline
		ABBCFHHH & $h_{1}=\{1, 18, 21, 42, 87, 106, 108, 119\}$, $h_{2}=\{1, 19, 27, 37, 78, 112,\
		
		113, 114\}$, $h_{3}=\{1, 20, 27, 32, 81, 113, 114, 115\}$\\
		\hline
		ABBDFGHH & $i_{1}=\{1, 18, 20, 51, 89, 101, 106, 108\}$, $i_{2}=\{1, 19, 21, 56, 85, 96,\
		
		106, 119\}$, $i_{3}=\{1, 19, 27, 46, 88, 96, 106, 112\}$, \\
		& $i_{4}=\{1, 20, 27, 46, 86, 101, 106,\
		
		115\}$\\
		\hline
		ABBEFFGH & $j_{1}=\{1, 18, 20, 64, 81, 89, 97, 113\}$, $j_{2}=\{1, 18, 20, 65, 83, 89, 95,\
		
		118\}$, $j_{3}=\{1, 18, 21, 64, 81, 87, 102, 116\}$, \\
		& $j_{4}=\{1, 18, 21, 65, 78, 87, 95, 111\}$,\
		
		$j_{5}=\{1, 19, 21, 64, 76, 85, 102, 109\}$, $j_{6}=\{1, 19, 21, 65, 78, 85, 100,\
		
		114\}$, \\
		& $j_{7}=\{1, 19, 27, 74, 76, 83, 102,\
		
		112\}$, $j_{8}=\{1, 20, 27, 70, 76, 83, 95, 115\}$\\
		\hline
		ABCCEFHH & $k_{1}=\{1, 18, 42, 43, 73, 83, 118, 119\}$, $k_{2}=\{1, 19, 34, 37, 64, 89,\
		
		109, 113\}$, $k_{3}=\{1, 20, 32, 35, 65, 85, 114, 118\}$, \\
		& $k_{4}=\{1, 21, 41, 42, 71, 76, 108,\
		
		109\}$\\
		\hline
		ABCDDGHH & $l_{1}=\{1, 18, 42, 50, 59, 97, 113, 119\}$, $l_{2}=\{1, 19, 37, 50, 56, 103,\
		
		113, 119\}$, $l_{3}=\{1, 20, 32, 51, 57, 94, 108, 114\}$, \\
		& $l_{4}=\{1, 21, 42, 48, 57, 100, 108,\
		
		114\}$, $l_{5}=\{1, 27, 42, 48, 59, 91, 113, 114\}$\\
		\hline
		ABCDEFGH & $m_{1}=\{1, 18, 43, 51, 73, 78, 101, 111\}$, $m_{2}=\{1, 19, 34, 51, 69, 89,\
		
		96, 106\}$, $m_{3}=\{1, 19, 43, 46, 64, 88, 94, 109\}$, \\
		& $m_{4}=\{1, 19, 43, 51, 69, 78, 94, 114\}$,\
		
		$m_{5}=\{1, 19, 43, 56, 74, 83, 94, 119\}$, $m_{6}=\{1, 20, 35, 56, 75, 85, 101, 106\}$,\
		
		\\
		& $m_{7}=\{1, 20, 41, 46, 65, 86, 103,\
		
		118\}$, $m_{8}=\{1, 20, 41, 51, 70, 76, 103, 108\}$, $m_{9}=\{1, 20, 41, 56, 75, 81,\
		
		103, 113\}$, \\
		& $m_{10}=\{1, 21, 41, 56, 71, 81, 96,\
		
		116\}$, $m_{11}=\{1, 27, 42, 46, 72, 87, 91, 106\}$\\
		\hline
		ABCEEFFG & $n_{1}=\{1, 18, 43, 64, 73, 81, 88, 92\}$, $n_{2}=\{1, 19, 34, 65, 74, 83, 89,\
		
		100\}$, $n_{3}=\{1, 20, 35, 64, 70, 76, 85, 97\}$, \\
		& $n_{4}=\{1, 21, 41, 65, 71, 78, 86, 105\}$,\
		
		$n_{5}=\{1, 27, 32, 71, 73, 81, 88, 96\}$, $n_{6}=\{1, 27, 32, 72, 74, 81, 87, 102\}$, \\
		& $n_{7}=\{1, 27, 37, 70, 72, 78, 87, 95\}$,\
		
		$n_{8}=\{1, 27, 37, 71, 73, 78, 86, 101\}$, $n_{9}=\{1, 27, 42, 71, 73, 76, 83, 91\}$\\
		\hline
		ABDDDGGH & $o_{1}=\{1, 18, 50, 51, 59, 101, 102, 116\}$, $o_{2}=\{1, 21, 48, 56, 57, 95,\
		
		96, 111\}$\\
		\hline
		ABDDEFGG & $p_{1}=\{1, 18, 50, 59, 65, 88, 92, 95\}$, $p_{2}=\{1, 19, 46, 50, 65, 88, 100,\
		
		103\}$, $p_{3}=\{1, 19, 50, 51, 69, 76, 102, 103\}$, \\
		& $p_{4}=\{1, 20, 46, 57, 64, 86, 94, 97\}$,\
		
		$p_{5}=\{1, 20, 56, 57, 75, 83, 94, 95\}$, $p_{6}=\{1, 21, 48, 57, 64, 86, 102, 105\}$,\
		
		\\
		& $p_{7}=\{1, 27, 48, 59, 70, 88, 95, 96\}$,\
		
		$p_{8}=\{1, 27, 48, 59, 74, 86, 101, 102\}$\\
		\hline
		ACCCDEHH & $q_{1}=\{1, 34, 41, 42, 59, 75, 109, 113\}$, $q_{2}=\{1, 35, 42, 43, 48, 69,\
		
		114, 118\}$, $q_{3}=\{1, 41, 42, 43, 46, 72, 109, 118\}$\\
		\hline
		ACCCEEFH & $r_{1}=\{1, 32, 34, 35, 65, 74, 87, 111\}$, $r_{2}=\{1, 34, 35, 37, 64, 70, 87,\
		
		116\}$, $r_{3}=\{1, 34, 35, 42, 69, 75, 87, 106\}$\\
		\hline
		ACCDDEGH & $s_{1}=\{1, 32, 34, 51, 57, 73, 96, 111\}$, $s_{2}=\{1, 34, 41, 51, 59, 70, 96,\
		
		116\}$, $s_{3}=\{1, 35, 37, 50, 56, 71, 101, 116\}$, \\
		& $s_{4}=\{1, 35, 43, 48, 56, 74, 101,\
		
		111\}$\\
		\hline
		ACCDEEFG & $t_{1}=\{1, 32, 35, 50, 65, 71, 88, 105\}$, $t_{2}=\{1, 34, 37, 57, 64, 73, 86,\
		
		92\}$, $t_{3}=\{1, 34, 41, 59, 65, 74, 86, 92\}$, \\
		& $t_{4}=\{1, 34, 42, 57, 73, 75, 83, 100\}$,\
		
		$t_{5}=\{1, 35, 42, 50, 69, 71, 76, 97\}$, $t_{6}=\{1, 35, 43, 48, 64, 70, 88, 105\}$, \\
		& $t_{7}=\{1, 41, 43, 51, 70, 72, 78, 105\}$,\
		
		$t_{8}=\{1, 41, 43, 56, 72, 74, 81, 92\}$\\
		\hline
		ACDDDEGG & $u_{1}=\{1, 32, 50, 51, 57, 72, 102, 105\}$, $u_{2}=\{1, 37, 50, 56, 57, 72,\
		
		92, 95\}$, $u_{3}=\{1, 42, 46, 50, 57, 72, 97, 100\}$\\
		\hline
		BBBBFFHH & $v_{1}=\{16, 19, 23, 27, 82, 88, 106, 117\}$\\
		\hline
		BBBCEFFH & $w_{1}=\{16, 19, 23, 43, 63, 80, 88, 109\}$, $w_{2}=\{16, 19, 27, 37, 62, 82,\
		
		90, 114\}$\\
		\hline
		BBBDDFGH & $x_{1}=\{16, 18, 20, 51, 52, 77, 99, 106\}$, $x_{2}=\{16, 18, 20, 52, 53, 79,\
		
		95, 118\}$, $x_{3}=\{16, 20, 24, 51, 57, 80, 99, 114\}$\\
		\hline
		BBBDEFFG & $y_{1}=\{16, 18, 20, 52, 63, 81, 90, 97\}$, $y_{2}=\{16, 19, 23, 50, 66, 79,\
		
		88, 100\}$, $y_{3}=\{16, 19, 27, 53, 74, 76, 82, 93\}$\\
		\hline
		BBCCDDHH & $z_{1}=\{16, 18, 36, 43, 49, 53, 118, 119\}$, $z_{2}=\{16, 19, 31, 34, 51, 52,\
		
		106, 117\}$, $z_{3}=\{16, 23, 41, 43, 49, 58, 109, 118\}$\\
		\hline
		BBCCDEFH & $a'_{1}=\{16, 18, 36, 43, 51, 62, 77, 111\}$, $a'_{2}=\{16, 19, 31, 43, 51, 62,\
		
		80, 114\}$, $a'_{3}=\{16, 19, 34, 37, 52, 63, 90, 109\}$, \\
		& $a'_{4}=\{16, 19, 37, 40, 50, 66, 90,\
		
		119\}$, $a'_{5}=\{16, 19, 40, 43, 53, 74, 80, 119\}$\\
		\hline
		BBCCEEFF & $b'_{1}=\{16, 18, 36, 43, 61, 63, 81, 88\}$, $b'_{2}=\{16, 20, 40, 41, 66, 75,\
		
		81, 90\}$\\
		\hline
		BBCDDEFG & $c'_{1}=\{16, 19, 31, 50, 51, 66, 76, 93\}$, $c'_{2}=\{16, 19, 34, 52, 53, 74,\
		
		79, 100\}$, $c'_{3}=\{16, 20, 40, 53, 57, 75, 80, 95\}$, \\
		& $c'_{4}=\{16, 20, 41, 47, 51, 66, 76,\
		
		99\}$\\
		\hline
		BBDDDDGG & $d'_{1}=\{16, 23, 49, 50, 57, 58, 97, 100\}$\\
		\hline
		BCCCDEEF & $e'_{1}=\{16, 34, 36, 37, 57, 61, 63, 86\}$, $e'_{2}=\{16, 40, 41, 43, 58, 61,\
		
		74, 81\}$\\
		\hline
		BCCDDDEG & $f'_{1}=\{16, 34, 36, 49, 53, 57, 75, 100\}$, $f'_{2}=\{16, 37, 40, 50, 57, 58,\
		
		61, 95\}$, $f'_{3}=\{16, 41, 43, 47, 51, 58, 62, 105\}$\\
		\hline
		CCCCDDEE & $g'_{1}=\{31, 32, 33, 34, 51, 52, 67, 73\}$\\
		\hline
	\end{tabular}
	\caption{Generators of Gosset's polytope. The generators are shown in the second column and their basis profiles in the first column. There are 135 generators but only 33 profiles, so many of the profiles have more than one generator associated with them. All generators with the same profile are labeled by the same letter and distinguished from each other by their subscripts. Each generator gives rise to fifteen bases and these bases, taken together, make up all the 2025 bases of Gosset's polytope.}
	\label{Tab 15} % is used to refer this table in the text
\end{table}

\begin{table}[ht]
	\centering % used for centering table
	\begin{tabular}{|c | c |} % centered columns (2 columns)
		\hline % inserts single horizontal line
		1-letter words (16) & $b_{1},e_{1},e_{2},g_{1},g_{2},g_{3},g_{4},g_{5},v_{1},z_{1},z_{2},z_{3},b'_{1},b'_{2},d'_{1},g'_{1}$ \\
		\hline
		2-letter words (22) & $a_{1}a_{2},a_{1}a_{3},a_{1}x_{1},a_{1}x_{2},a_{1}x_{3},c_{1}c_{2},$\\
		& $d_{1}c'_{1},d_{1}c'_{2},d_{1}c'_{3},d_{1}c'_{4},d_{1}d_{2},d_{1}d_{3},$\\
		& $d_{1}d_{4},d_{1}d_{5},d_{1}d_{6},d_{1}d_{7},d_{1}d_{8},h_{1}h_{2},$\\
    	& $h_{1}h_{3},h_{1}r_{1},h_{1}r_{2},h_{1}r_{3}$\\ 
		\hline
		3-letter words (56) & $a_{1}c_{1}w_{1},a_{1}c_{1}w_{2},a_{1}c_{1}e'_{1},a_{1}c_{1}e'_{2},a_{1}h_{1}l_{1},a_{1}h_{1}l_{2},$\\
		& $a_{1}h_{1}l_{3},a_{1}h_{1}l_{4},a_{1}h_{1}l_{5},a_{1}h_{1}n_{1},a_{1}h_{1}n_{2},a_{1}h_{1}n_{3},$\\
		& $a_{1}h_{1}n_{4},a_{1}h_{1}n_{5},a_{1}h_{1}n_{6},a_{1}h_{1}n_{7},a_{1}h_{1}n_{8},a_{1}h_{1}n_{9},$\\
		& $c_{1}d_{1}f_{1},c_{1}d_{1}f_{2},c_{1}d_{1}f_{3},c_{1}d_{1}f_{4},c_{1}d_{1}f_{5},c_{1}d_{1}f_{6},$\\
		& $c_{1}d_{1}f_{7},c_{1}d_{1}f_{8},c_{1}d_{1}f_{9},c_{1}d_{1}a'_{1},c_{1}d_{1}a'_{2},c_{1}d_{1}a'_{3},$\\
		& $c_{1}d_{1}a'_{4},c_{1}d_{1}a'_{5},c_{1}h_{1}i_{1},c_{1}h_{1}i_{2},c_{1}h_{1}i_{3},c_{1}h_{1}i_{4},$\\ 
		& $c_{1}h_{1}t_{1},c_{1}h_{1}t_{2},c_{1}h_{1}t_{3},c_{1}h_{1}t_{4},c_{1}h_{1}t_{5},c_{1}h_{1}t_{6},$\\
		& $c_{1}h_{1}t_{7},c_{1}h_{1}t_{8},d_{1}h_{1}j_{1},d_{1}h_{1}j_{2},d_{1}h_{1}j_{3},d_{1}h_{1}j_{4},$\\
		& $d_{1}h_{1}j_{5},d_{1}h_{1}j_{6},d_{1}h_{1}j_{7},d_{1}h_{1}j_{8},d_{1}h_{1}s_{1},d_{1}h_{1}s_{2},$\\
		& $d_{1}h_{1}s_{3},d_{1}h_{1}s_{4}$\\
		\hline
		4-letter words (26) & $a_{1}c_{1}d_{1}f'_{1},a_{1}c_{1}d_{1}f'_{2},a_{1}c_{1}d_{1}f'_{3},a_{1}c_{1}d_{1}y_{1},$\\
		& $a_{1}c_{1}d_{1}y_{2},a_{1}c_{1}d_{1}y_{3},a_{1}c_{1}h_{1}o_{1},a_{1}c_{1}h_{1}o_{2}, $\\
		& $a_{1}d_{1}h_{1}k_{1},a_{1}d_{1}h_{1}k_{2},a_{1}d_{1}h_{1}k_{3},a_{1}d_{1}h_{1}k_{4}, $\\
		& $a_{1}d_{1}h_{1}p_{1},a_{1}d_{1}h_{1}p_{2},a_{1}d_{1}h_{1}p_{3},a_{1}d_{1}h_{1}p_{4}, $\\ 
		& $a_{1}d_{1}h_{1}p_{5},a_{1}d_{1}h_{1}p_{6},a_{1}d_{1}h_{1}p_{7},a_{1}d_{1}h_{1}p_{8}, $\\ 
		& $c_{1}d_{1}h_{1}q_{1},c_{1}d_{1}h_{1}q_{2},c_{1}d_{1}h_{1}q_{3},c_{1}d_{1}h_{1}u_{1}, $\\
		& $c_{1}d_{1}h_{1}u_{2},c_{1}d_{1}h_{1}u_{3} $\\
		\hline
		5-letter words (11) & $a_{1}c_{1}d_{1}h_{1}m_{1},a_{1}c_{1}d_{1}h_{1}m_{2},a_{1}c_{1}d_{1}h_{1}m_{3}, $\\ 
		& $a_{1}c_{1}d_{1}h_{1}m_{4},a_{1}c_{1}d_{1}h_{1}m_{5},a_{1}c_{1}d_{1}h_{1}m_{6}, $\\
		& $a_{1}c_{1}d_{1}h_{1}m_{7},a_{1}c_{1}d_{1}h_{1}m_{8},a_{1}c_{1}d_{1}h_{1}m_{9}, $\\
		& $a_{1}c_{1}d_{1}h_{1}m_{10},a_{1}c_{1}d_{1}h_{1}m_{11} $\\
		\hline
	\end{tabular}
	\caption{The 131 words representing a linearly independent set of vectors in the nullspace of the matrix $M^{T}$ of Gosset's polytope. The 16 single letter words all give $30_{2}-15_{4}$ parity proofs and the 3- and 5-letter words represent parity proofs of 45 and 75 bases, respectively. All the parity proofs in this table are both minimal and irreducible.}
	\label{Tab9} % is used to refer this table in the text
\end{table}

%% \begin{table}[ht]
%% 	\centering % used for centering table
%% 	\begin{tabular}{|c | c | c| c |} % centered columns (2 columns)
%% 		\hline % inserts single horizontal line
%% 		\bf{Length} & \bf{Proofs} & \bf{Length} & \bf{Proofs }\\
%% 		\hline
%% 		& & & \\
%% 		\hline
%% 		& & & \\
%% 		\hline
%% 		& & & \\
%% 		\hline
%% 		& & & \\
%% 		\hline
%% 		& & & \\
%% 		\hline
%% 		& & & \\
%% 		\hline
%% 		& & & \\
%% 		\hline
%% 		& & & \\
%% 		\hline
%% 		& & & \\
%% 		\hline
%% 		& & & \\
%% 		\hline
%% 	\end{tabular}
%% 	\caption{Fifteen-fold symmetric parity proofs of Gosset's polytope. Shown above are the number of proofs for each word of odd length from 1 to 39.}
%% 	\label{Tab11} % is used to refer this table in the text
%% \end{table}

\begin{table}[ht]
\begin{tabular}{c c}
\end{tabular}

\end{table}

\begin{table}[th]
  \begin{tabular}{|c |c |c |c|}
    \hline
    \textbf{Length} & \textbf{Proofs} & \textbf{Length} & \textbf{Proofs}\\
    \hline
    1 &  16 & 69 &  180520872122823098245862570199127030680 \\
    \hline
    3 &  25812 & 71 &  155821839317285933878889210179604711024 \\
    \hline
    5 &  21653868 & 73 &  119534561668054963002852653943661854832 \\
    \hline
    7 &  8652003024 & 75 &  81455804005150246833959629030543329328 \\
    \hline
    9 &  1953439358160 & 77 &  49274358540367715978978010339417431984 \\
    \hline
    11 &  279696462166032 & 79 &  26436389051355999506519101553074086512 \\
    \hline
    13 &  27345712822595216 & 81 &  12565444178730938040258584478235056624 \\
    \hline
    15 &  1922273370347815632 & 83 &  5283911436898023021008844055037541520 \\
    \hline
    17 &  100919351992599919440 & 85 &  1962595676562122833284741838569192272 \\
    \hline
    19 &  4073954893763530856880 & 87 &  642656964391499727733011021261899760 \\
    \hline
    21 &  129396567337947795186928 & 89 &  185116714972832401789642506963347440 \\
    \hline
    23 &  3294242253859564715609424 & 91 &  46787741146979619043486256387889200 \\
    \hline
    25 &  68256699499928681794401616 & 93 &  10346237289631303832281044955692720 \\
    \hline
    27 &  1165808870376411660238444048 & 95 &  1995097493028567348842572411414000 \\
    \hline
    29 &  16591240524716074401300742544 & 97 &  334230250120765120581947447514800 \\
    \hline
    31 &  198559684989345477675990749008 & 99 &  48436171064707869561000746697120 \\
    \hline
    33 &  2014177410611613275385516844560 & 101 &  6042532231834847059425191542336 \\
    \hline
    35 &  17437021583294837611603671701504 & 103 &  645318976215369969425403092912 \\
    \hline
    37 &  129599484740704840232429138427040 & 105 &  58622383187330985956893984176 \\
    \hline
    39 &  831290622095236372451283545241360 & 107 &  4496691357165952988389892464 \\
    \hline
    41 &  4622786874090582681904187711459600 & 109 &  288778344038199652858660080 \\
    \hline
    43 &  22376745765946773959473914330329680 & 111 &  15373130518005207330609072 \\
    \hline
    45 &  94615209875003228427091933876846800 & 113 &  670509485297303171419696 \\
    \hline
    47 &  350540162395363486324657008310384400 & 115 &  23628938383825022644112 \\
    \hline
    49 &  1141043998001234207494166219135912080 & 117 &  661582418589474325200 \\
    \hline
    51 &  3270992794270204725504176218418242608 & 119 &  14417050288814344240 \\
    \hline
    53 &  8274804703066715294324920218378861296 & 121 &  238298351913548848 \\
    \hline
    55 &  18505472335949199660124342087904997008 & 123 &  2890197197138544 \\
    \hline
    57 &  36639907632581122125005175119122150032 & 125 &  24613298133232 \\
    \hline
    59 &  64307213688276510679781711925729104720 & 127 &  138432535088 \\
    \hline
    61 &  100150578694856860867732458243810022096 & 129 &  469443904 \\
    \hline
    63 &  138508301615365274584106443908685446288 & 131 &  829828 \\
    \hline
    65 &  170205393715804635473398600245423862256 & 133 &  540 \\
    \hline
    67 &  185909554872758116099101213433344650760 &  &  \\
    \hline
  \end{tabular}
  \caption{Fifteen-fold symmetric parity proofs of Gosset's polytope. Shown above are the number of proofs for each word of odd length from 1 to 133. The total number of proofs is $2^{130}$.}
\end{table}